\documentclass[12pt, twoside, fleqn, a4paper]{article}

\usepackage{latexsym}  
\usepackage{amscd}    
\usepackage{theorem}  
\usepackage{pifont}  
\usepackage{mathbbol} 
\usepackage{amsfonts}  
\usepackage{xspace}  
\usepackage{amssymb} 
\usepackage{fancyhdr} 
\usepackage{mathrsfs} 
\usepackage{amsmath} 
\usepackage{graphics} 
\usepackage{graphicx} 
\usepackage{euscript}
\usepackage{psfrag}
\usepackage{empheq} 
\usepackage{mathabx} 
\usepackage[parfill]{parskip}     
\usepackage[colorlinks=true, allcolors=blue]{hyperref}
\usepackage{cancel} 
\usepackage{bigints} 

\title{Large deviation of preimages of points under the dynamics of a holomorphic correspondence}
\author{Bharath Krishna Seshadri\footnote{Indian Institute of Science Education and Research Thiruvananthapuram (IISER-TVM). \\ ORCID: 0009 - 0004 - 0541 - 8040\ \ email: \texttt{bharathmaths21@iisertvm.ac.in}}\ \ and\ \ Shrihari Sridharan\footnote{Indian Institute of Science Education and Research Thiruvananthapuram (IISER-TVM). \\ ORCID: 0000-0003-2434-4767\ \ email: \texttt{shrihari@iisertvm.ac.in}\ (Corresponding Author) \\ The second named author thanks the support provided by NBHM through a research grant \\ No. 02011/35/2025/NBHM(R.P)/R\&D II/9832}}

\DeclareFontFamily{OT1}{pzc}{}
\DeclareFontShape{OT1}{pzc}{m}{it}%
              {<-> s * [0.900] pzcmi7t}{}
\DeclareMathAlphabet{\mathpzc}{OT1}{pzc}%
                                 {m}{it}

{\theorembodyfont{\slshape} \newtheorem{theorem}{Theorem}[section]}
{\theorembodyfont{\slshape} \newtheorem{definition}[theorem]{Definition}}
{\theorembodyfont{\slshape} \newtheorem{lemma}[theorem]{Lemma}}
{\theorembodyfont{\slshape} }
{\theorembodyfont{\slshape} }
{\theorembodyfont{\slshape} \newtheorem{corollary}[theorem]{Corollary}} 
{\theorembodyfont{\slshape} } 
{\theorembodyfont{\slshape} }
{\theorembodyfont{\slshape} }

\numberwithin{equation}{section}

\newenvironment{proof}{\paragraph{Proof:}}{\hfill$\bullet$}

\begin{document}

\maketitle

\begin{abstract} 
\noindent 
Our main objective in this paper is to prove a large deviation result for the distribution of preimages of a point under the iteration of a holomorphic correspondence with respect to the unique equilibrium state associated to a H\"{o}lder continuous potential defined on the support of the Dinh-Sibony measure of the correspondence. As a consequence, we obtain a convergence result which generalises the result on equidistribution of preimages by Dinh and Sibony. We also study a situation where the Dinh-Sibony measure happens to be the unique measure of maximal entropy for an appropriate holomorphic correspondence and also record some comments that connect this work with some recent works in the literature.
\end{abstract}

\begin{tabular}{l l} 
{\bf Keywords} & Large deviations \\ 
& Unique equilibrium state \\ 
& Distance expanding holomorphic correspondence \\
& Dinh-Sibony measure and its support \\ 
& \\
{\bf MSC Subject} & 37A35, 37A50, 37F05, 37F15, 37F80 \\ 
{\bf Classifications} & 
\end{tabular} 
\bigskip

\section{Introduction} 
\label{sec:one}

The study of dynamics of holomorphic correspondences is now a rapidly evolving branch of dynamics that touches upon complex dynamics, geometry, and ergodic theory. The current work focuses on the complex analytic and mostly ergodic theoretic perspectives by studying large deviation results on the distribution of preimages of a point under the dynamics of a holomorphic correspondence. We shall work with holomorphic correspondences on $\widehat{\mathbb{C}}$, not merely because things are nice and convenient on $\widehat{\mathbb{C}}$, but also because it is much more enlightening, in the sense that holomorphic correspondences on $\widehat{\mathbb{C}}$ are generalisations of rational maps on $\widehat{\mathbb{C}}$, thus providing a nice analogous set up to fall back to in order to understand things better. 

There have been many works related to large deviation theory in complex dynamics. In \cite{lop:1990}, Lopes investigated the connection between entropy and large deviation. This was one of the earliest works to study large deviations using thermodynamic formalism in the case of rational maps. Lopes' work focused on studying hyperbolic rational maps and the measure of maximal entropy associated with them. Large deviation principles for non-uniformly hyperbolic rational maps was studied in \cite{cr:2011}, where the authors proved some large deviation results by investigating thermodynamic formalism for Topological Collet-Eckmann maps. The works \cite{ps:1996} and \cite{ps:2007} respectively dealt with the study of large deviation results with regards to the distribution of preimages and periodic points under the dynamics of a rational map. In \cite{ssa:2021} the distribution of preimages and periodic points with respect to the dynamics of the skew product map associated with a rational semigroup was studied. Our work is mainly inspired by the ideas in \cite{ps:1996} and hence we now state the main result and its corollary found therein, before proceeding to brief about the analogous set up in the setting of holomorphic correspondences, that we study in our work. 

In the following theorem, $J_{R}$ denotes the Julia set of the rational map $R$ defined on $\widehat{\mathbb{C}}$, $\mathcal{M}$ the space of Borel probability measures on $J_{R}$ and ${\rm Pr} (R, f)$ denotes the pressure of the potential $f$ with respect to $R$. 
\begin{theorem}\cite{ps:1996}
\label{thm:psld}
    Let $f \in \mathcal{C} (J_{R}, \mathbb{R})$ be H\"{o}lder continuous with ${\rm Pr} (R, f) \, > \, \sup f$ and let $\mu$ be the unique equilibrium state for $f$. Let $x \in J_{R}$. Then for any weak* open neighbourhood $\mathcal{U} \subset \mathcal{M}$ of $\mu$ we have that the weighted proportion of the measures $\mu_{y, \, n} = \frac{1}{n} \left( \delta_{y} + \delta_{Ry} + \cdots + \delta_{R^{n - 1}y} \right) \, \notin \, \mathcal{U}$ tends to zero exponentially fast, in the sense that, there exists $C > 0$ and $0 < \eta < 1$ such that,
    \[
    \dfrac{\displaystyle \sum \limits_{\begin{minipage}{1.5cm} \begin{center} \scriptsize{$R^{n}y \, = \, x$} \\ \scriptsize{$\mu_{y, \, n} \, \notin \, \mathcal{U}$} \end{center} \end{minipage}} e^{f^{n}(y)}}{\displaystyle \sum \limits_{R^{n}y \, = \, x} e^{f^{n}(y)}} \, \le \, C\eta^{n}, \ \ \ \text{for} \ n \ge 0,
    \]
    where we denote $f^{n}(y) = f(y) + f(Ry) + \cdots + f(R^{n - 1}y)$.
\end{theorem}

Pollicott and Sharp obtained two results as corollaries to the above theorem (see Corollaries 1 and 2 in \cite{ps:1996}) that generalised the result on equidistribution of preimages by Lyubich.

Now, turning our attention to holomorphic correspondences, we start with one of the breakthrough results in this study, that came courtesy the work \cite{ds:2006} by Dinh and Sibony. Although the theorem holds true in a much more general context, we state a version of it that suits our needs in this paper. The basics and notations involved in the following theorem shall be elaborated in Section \ref{sec:two}.

\begin{theorem}\cite{ds:2006}
\label{thm:ds} 
Let $\digamma$ be a holomorphic correspondence on $\widehat{\mathbb{C}}$ with $d_{{\rm top}} > d_{{\rm fwd}}$. Then, there is an exceptional set $\mathcal{E} \subsetneq \widehat{\mathbb{C}}$ such that 
\[ \frac{1}{d_{{\rm top}}^{n}} \left( \digamma^{n} \right)^{*} \delta_{z_{0}}\ \ \rightharpoonup\ \ \mu_{\digamma},\ \ \ \ \text{in the weak* topology,\ for all}\ z_{0} \in \widehat{\mathbb{C}} \setminus \mathcal{E}. \] 
Here, $\left( \digamma^{n} \right)^{*} \delta_{z_{0}}$ denotes the pullback of the Dirac delta measure $\delta_{z_{0}}$ under the holomorphic correspondence $\digamma^{n}$. Moreover, $\mu_{\digamma}$ satisfies $\digamma^{*} \mu_{\digamma} = d_{{\rm top}} \, \mu_{\digamma}$.
\end{theorem} 

The convergence in Theorem \ref{thm:ds} can be rewritten in the following way using the notations in Section \ref{sec:two}.

\begin{equation}
\label{eqn:dscgt}
    \frac{1}{d_{{\rm top}}^{n}} \sum_{\mathfrak{X}^{-}_{n} \left( z_{0}; \boldsymbol{\beta} \right)_{\boldsymbol{k}}\, \in\, \mathscr{Q}_{n}^{\Gamma} \left( \widehat{\mathbb{C}} \right)} \delta_{\Pi_{(n,n)}^{-} \left( \mathfrak{X}^{-}_{n} \left( z_{0}; \boldsymbol{\beta} \right)_{\boldsymbol{k}} \right)} \ \ \rightharpoonup\ \ \mu_{\digamma}, \ \ \ \forall z_{0} \, \in \, \widehat{\mathbb{C}} \setminus \mathcal{E}.
\end{equation}

The support of the measure $\mu_{\digamma}$ from Theorem \ref{thm:ds} is known to contain some dynamically interesting properties (as one may find the details in \cite{bs:2016}) and is the analogue of the Julia set from the study of rational maps. We denote this set by $\Omega_{{\rm DS}}$ and call it the \emph{support of the Dinh-Sibony measure} of $\digamma$. In this work, we will mostly be working with correspondences for which such a set $\Omega_{{\rm DS}}$ exists, as we will be restricting the dynamics of $\digamma$ to this space. This restriction is possible due to a result by Londhe (see Theorem 1.8 in \cite{lon:2022}) where the author proved that $\Omega_{{\rm DS}}$ is backward invariant and that there is at least one forward orbit of any point in $\Omega_{{\rm DS}}$ that is fully inside $\Omega_{{\rm DS}}$. Results mirroring Theorem \ref{thm:ds} have been obtained for some other classes of holomorphic correspondences in various works such as \cite{bs:2016, dkw:2020, hem:2024, vd:2024} and hence it is possible to define $\Omega_{{\rm DS}}$ for such correspondences too.

The large deviation problem also requires a considerable amount of thermodynamic formalism. This particular aspect in the dynamics of holomorphic correspondences was first studied in detail in \cite{ds:2008}. In that work, Dinh and Sibony also showed that the topological entropy of a holomorphic correspondence is finite, which is an essential characteristic for our problem. The concept of measure theoretic entropy and pressure were later developed in \cite{ss:2025, ss:2026}. One crucial element in Theorem \ref{thm:psld} is the existence of a unique equilibrium state for a given H\"{o}lder continuous potential on $\Omega_{{\rm DS}}$. One cannot expect this to exist under the dynamics of any arbitrary holomorphic correspondence, but in \cite{bks:2026}, the authors identified a class of holomorphic correspondences under the dynamics of which, each H\"{o}lder continuous potential has a unique equilibrium state. Thus, the stage is set for the work on large deviations and hence, the study in this paper can be seen as an application of the theory developed in \cite{bks:2026} as well. We now state the main results of our work. The relevant notations are explained in detail in Sections \ref{sec:two} and \ref{sec:three}.

For a point $z_0 \in \Omega_{{\rm DS}}$ and $\mathfrak{X}_{n}^{-} \left( z_{0}; \boldsymbol{\beta} \right)_{\boldsymbol{k}} \in \mathscr{Q}_{n}^{\Gamma} (z_0)$, we set,
\begin{equation}
\label{eqn:muzee}
    \mu_{\mathfrak{X}_{n}^{-} \left( z_{0}; \boldsymbol{\beta} \right)_{\boldsymbol{k}}} = \frac{1}{n} \sum_{r \, =\, 1}^{n} \delta_{\Pi_{(r, n)}^{-} \left( \mathfrak{X}_{n}^{-} \left( z_{0}; \boldsymbol{\beta} \right)_{\boldsymbol{k}} \right)} \ \ \text{and} \ \  \mathcal{Z}_{n}^{-} (f) \left( \mathfrak{X}_{n}^{-} \left( z_{0}; \boldsymbol{\beta} \right)_{\boldsymbol{k}} \right) = \sum \limits_{r \, = \, 1}^{n} f\left(\Pi_{(r, n)}^{-} \left( \mathfrak{X}_{n}^{-} \left( z_{0}; \boldsymbol{\beta} \right)_{\boldsymbol{k}} \right)\right).
\end{equation}

\begin{theorem} 
\label{thm:main} 
Let $\digamma$ be a holomorphic correspondence which is distance expanding on $\Omega_{{\rm DS}}$ and let $f \in \mathcal{C}(\Omega_{{\rm DS}}, \mathbb{R})$ be a H\"{o}lder continuous function satisfying the condition ${\rm Pr} \left( \Gamma,\; f \right) \, > \, \sup f$, with its unique equilibrium state being denoted by $\nu_{f} \in \mathscr{S}^{\Gamma}$. Then, for any $z_{0} \in \Omega_{{\rm DS}}$ and for any weak*-open neighbourhood $\mathcal{V} \subsetneq \mathcal{M} \left( \Omega_{{\rm DS}} \right)$ of $\nu_{f}$, there exist positive real constants $K$ and $r$ with $0 < r < 1$ such that 
\[ \dfrac{\displaystyle \sum\limits_{\begin{minipage}{3cm} \begin{center} \scriptsize{$ \mathfrak{X}_{n}^{-} \left( z_{0}; \boldsymbol{\beta} \right)_{\boldsymbol{k}} \, \in\,  \mathscr{Q}_{n}^{\Gamma} \left( z_{0} \right)$} \\ \scriptsize{$\mu_{\mathfrak{X}_{n}^{-} \left( z_{0}; \boldsymbol{\beta} \right)_{\boldsymbol{k}}} \, \notin \,  \mathcal{V}$} \end{center} \end{minipage}} e^{\mathcal{Z}_{n}^{-} (f) \left( \mathfrak{X}_{n}^{-} \left( z_{0}; \boldsymbol{\beta} \right)_{\boldsymbol{k}} \right)}}{\displaystyle \sum\limits_{\mathfrak{X}_{n}^{-} \left( z_{0}; \boldsymbol{\beta} \right)_{\boldsymbol{k}}\, \in\, \mathscr{Q}_{n}^{\Gamma} \left( z_{0} \right)} e^{\mathcal{Z}_{n}^{-} (f) \left( \mathfrak{X}_{n}^{-} \left( z_{0}; \boldsymbol{\beta} \right)_{\boldsymbol{k}} \right)}}\ \ \le\ \ K r^{n}, \ \ \ \ \forall n \ge 0. \] 
\end{theorem} 

Taking $f \equiv 0$ in Theorem \ref{thm:main}, we obtain the following corollary. 

\begin{corollary} 
Let $z_{0} \in \Omega_{{\rm DS}}$ and $\mu_{0} \in \mathscr{S}^{\Gamma}$ denote the unique measure of maximal entropy of holomorphic correspondence $\digamma$, which is distance expanding on $\Omega_{{\rm DS}}$. Then for any weak*-open neighbourhood $\mathcal{V} \subsetneq \mathcal{M} \left( \Omega_{{\rm DS}} \right)$ of $\mu_{0}$, there exist positive real constants $K$ and $r$ with $0 < r < 1$ such that 
\[ \frac{1}{d_{{\rm top}}^{n}} \# \left\{ \mathfrak{X}_{n}^{-} \left( z_{0}; \boldsymbol{\beta} \right)_{\boldsymbol{k}} \in \mathscr{Q}_{n}^{\Gamma} \left( z_{0} \right)\ :\ \mu_{\mathfrak{X}_{n}^{-} \left( z_{0}; \boldsymbol{\beta} \right)_{\boldsymbol{k}}} \notin \mathcal{V} \right\}\ \ \le\ \ K r^{n},\ \ \ \forall n \ge 0. \] 
\end{corollary} 

Another interesting consequence of Theorem \ref{thm:main} is the following theorem, which in some sense generalises the expression in Equation \eqref{eqn:dscgt} by considering a weighted proportion involving measures. The limit measure in this result turns out to be the unique equilibrium state associated to the potential involved. 

\begin{theorem} 
\label{thm:muf}
Let $\digamma$ be a holomorphic correspondence which is distance expanding on $\Omega_{{\rm DS}}$ and let $f \in \mathcal{C}(\Omega_{{\rm DS}}, \mathbb{R})$ be a H\"{o}lder continuous function satisfying ${\rm Pr} (\Gamma, f) > \sup f$, whose unique equilibrium state is denoted by $\nu_{f} \in \mathscr{S}^{\Gamma}$. Let $z_{0} \in \Omega_{{\rm DS}}$. Then, 
\[ \dfrac{\sum\limits_{\mathfrak{X}_{n}^{-} \left( z_{0}; \boldsymbol{\beta} \right)_{\boldsymbol{k}}\, \in\, \mathscr{Q}_{n}^{\Gamma} \left( z_{0} \right)} e^{\mathcal{Z}_{n}^{-} (f) \left( \mathfrak{X}_{n}^{-} \left( z_{0}; \boldsymbol{\beta} \right)_{\boldsymbol{k}} \right)} \mu_{\mathfrak{X}_{n}^{-} \left( z_{0}; \boldsymbol{\beta} \right)_{\boldsymbol{k}}}}{\sum\limits_{\mathfrak{X}_{n}^{-} \left( z_{0}; \boldsymbol{\beta} \right)_{\boldsymbol{k}}\, \in\, \mathscr{Q}_{n}^{\Gamma} \left( z_{0} \right)} e^{\mathcal{Z}_{n}^{-} (f) \left( \mathfrak{X}_{n}^{-} \left( z_{0}; \boldsymbol{\beta} \right)_{\boldsymbol{k}} \right)}}\ \ \to\ \ \nu_{f},\ \ \ \text{in the weak*-topology, as}\ n \to \infty. \] 
\end{theorem} 

The above theorem has the following interesting corollary which shall be proved in Section \ref{sec:six}. Theorem \ref{thm:muf} also acts as a nice generalisation to the result regarding the equidistribution of preimages of points by Dinh and Sibony, which was proved in \cite{ds:2006}. A few more details on this result and about the Dinh-Sibony measure itself will be discussed in Section \ref{sec:six}.

\begin{corollary}
\label{cor:dsisuniqueent}
    Let $\digamma$ be a holomorphic correspondence which is distance expanding on $\Omega_{{\rm DS}}$. Then, the Dinh-Sibony measure of $\digamma$ is the unique measure of maximal entropy for $\digamma$.
\end{corollary}

The paper is organised as follows. In Section \ref{sec:two}, we state the preliminaries on the dynamics of holomorphic correspondences on $\widehat{\mathbb{C}}$. We mainly focus on how its iterates are defined and we also study the space of orbits of points in $\widehat{\mathbb{C}}$ under the dynamics of the given correspondence. In Section \ref{sec:three}, we state the results on ergodic theory and thermodynamic formalism for correspondences that are necessary for our work. This paves the way for us to relate the dynamics of the correspondence and the associated shift map $\sigma^{\Gamma}$, in Section \ref{sec:four}. In that section, we prove more results on thermodynamic formalism that are essential for the proof of our main results. Section \ref{sec:five} is fully dedicated to state the proofs of our main results, namely Theorems \ref{thm:main} and \ref{thm:muf}. In Section \ref{sec:six}, we prove Corollary \ref{cor:dsisuniqueent} and explore some properties of the Dinh-Sibony measure and mention some comments that connect our work with some other works in the field. We end this manuscript with Section \ref{sec:seven}, which contains a short discussion on the strategy used to prove Lemma \ref{lem:lf}.

\section{Holomorphic correspondences and their dynamics}
\label{sec:two}

In this section, we introduce the basic setting on which this paper is built. Most of the material mentioned in this section are well established in the literature by now and we shall state the required results briefly. One may refer \cite{ds:2006, ds:2008, bs:2016, lon:2024, ss:2025, bks:2026, ss:2026} for more details on the concepts we introduce here in this section. 

Let $\widehat{\mathbb{C}}$ denote the Riemann sphere, which is a compact metric space with respect to the spherical metric $\rho$. Denote by $\pi_{1}$ and $\pi_{2}$ the canonical projections which maps each point of $\widehat{\mathbb{C}} \times \widehat{\mathbb{C}}$ onto their first and second coordinates respectively. A \emph{holomorphic $1$-chain} in $\widehat{\mathbb{C}} \times \widehat{\mathbb{C}}$ is a formal sum of finitely many (not necessarily distinct) irreducible complex subvarieties of $\widehat{\mathbb{C}} \times \widehat{\mathbb{C}}$ with each of them having dimension $1$. A holomorphic $1$-chain  $\Gamma = \sideset{}{'} \sum \limits_{j \, = \, 1}^{N} \Gamma_{j}$ is said to be a \emph{holomorphic correspondence} if for each $1 \le j \le N$, we have $\left. \pi_{1} \right|_{\Gamma_{j}}$ and $\left. \pi_{2} \right|_{\Gamma_{j}}$ to be surjective and the sets $\bigcup\limits_{j\, =\, 1}^{N} \pi_{2} \left( \pi_{1}^{-1} \{z_{0}\} \cap \Gamma_{j} \right)$ and $\bigcup\limits_{j\,=\, 1}^{N} \pi_{1} \left( \pi_{2}^{-1} \{z_{0}\} \cap \Gamma_{j} \right)$ to be finite for any $z_{0} \in \widehat{\mathbb{C}}$. The expression denoted by $\Gamma$ above is said to be the \emph{variety-representation} of the holomorphic correspondence and the notation $\sideset{}{'} \sum$ denotes that the varieties in the formal sum are allowed to repeat, according to their multiplicities. The correspondence $\Gamma$ gives rise to a set-valued map, which we shall denote by $\digamma$. For any $A \subset \widehat{\mathbb{C}}$, we define $\digamma (A) = \bigcup\limits_{j\, =\, 1}^{N} \pi_{2} \left( \pi_{1}^{-1} (A) \cap \Gamma_{j} \right)$ and $\digamma^{\dagger} (A) = \bigcup\limits_{j\,=\, 1}^{N} \pi_{1} \left( \pi_{2}^{-1} (A) \cap \Gamma_{j} \right)$. It must be noted that these `sets' allow elements to repeat according to multiplicity. The set valued map $\digamma^{\dagger}$ defines the holomorphic correspondence $\Gamma^{\dagger} = \sideset{}{'} \sum \limits_{j \, = \, 1}^{N} \Gamma_{j}^{\dagger}$, where $\Gamma_j^{\dagger} = \{(z,w) \in \widehat{\mathbb{C}} \times \widehat{\mathbb{C}} : (w, z) \in \Gamma_j \}$.  The cardinality of $\digamma (z_0)$ and that of $\digamma^{\dagger} (z_0)$ for a generic point $z_0 \in \widehat{\mathbb{C}}$ is called the \emph{forward} and \emph{topological degree} respectively for the holomorphic correspondence $\Gamma$; the former being denoted by $d_{{\rm fwd}}$ while the latter being denoted by $d_{{\rm top}}$. The sets $\digamma^{n} (A)$ and $(\digamma^{\dagger})^n (A)$ can be defined recursively and more details on the basics of their iterative dynamics can be found in \cite{bs:2016}. 

There is another way to study the iterative dynamics of $\Gamma$, which is done by studying the space of orbits of points in $\widehat{\mathbb{C}}$ under the iterates of $\Gamma$. One of the earliest works to adopt this perspective was \cite{ds:2008} and the notations that we use in this work are akin to the ones in \cite{ss:2026, bks:2026}. For a given $n \in \mathbb{Z}_{+}$, we begin by defining the spaces of backward orbits of length $n$ under the dynamics of $\Gamma$. 
\begin{eqnarray*} 
\mathscr{Q}_{n}^{\Gamma} \left( z_{0} \right) & = & \Big\{ \left( z_{k_{n}}^{(-n)}, \cdots, z_{k_{1}}^{(-1)}, z_{0};\; \beta_{n}, \cdots, \beta_{1} \right)\ :\ \left(z_{k_{r}}^{(-r)}, z_{k_{r - 1}}^{(-(r - 1))} \right) \in \Gamma_{\beta_{r}} \ \text{where}\ z_{k_{0}}^{(0)} = z_{0}, \\ 
& & \hspace{+0.8cm} \beta_{r} \in \left\{ 1, 2, \cdots, N \right\}\ \text{and}\ 1 \le k_{r} \le \delta_{\beta_{r}} \left( z_{k_{r - 1}}^{(-(r - 1))} \right)\ \text{for}\ n \ge r \ge 1 \Big\},\ \ \text{and} \\ 
\mathscr{Q}_{n}^{\Gamma} \left( \widehat{\mathbb{C}} \right) & = & \bigcup_{z_{0}\, \in\, \widehat{\mathbb{C}}} \mathscr{Q}_{n}^{\Gamma} \left( z_{0} \right). 
\end{eqnarray*} 
The symbol $\delta_{\beta_{r}} \left( z_{k_{r - 1}}^{(-(r - 1))} \right)$ denotes the cardinality of the set $\digamma^{\dagger} \left(z_{k_{r - 1}}^{(-(r - 1))} \right) \cap \pi_{1} (\Gamma_{\beta_r})$. We represent any point in $\mathscr{Q}_{n}^{\Gamma} \left( z_{0} \right)$ by $\mathfrak{X}_{n}^{-} \left( z_{0}; \boldsymbol{\beta} \right)_{\boldsymbol{k}}$ where $\boldsymbol{\beta} = \left( \beta_{n}, \cdots, \beta_{1} \right)$ and $\boldsymbol{k} = \left( k_{n}, \cdots, k_{1} \right)$ and $z_0 \in \widehat{\mathbb{C}}$. Moreover, let 
\begin{eqnarray*} 
\mathscr{Q}^{\Gamma} \left( \widehat{\mathbb{C}} \right) & = & \Big\{ \left( \cdots, z_{k_{2}}^{(-2)}, z_{k_{1}}^{(-1)}, z_{0};\; \cdots, \beta_{2}, \beta_{1} \right)\ :\ \left(z_{k_{r}}^{(-r)}, z_{k_{r - 1}}^{(-(r - 1))} \right) \in \Gamma_{\beta_{r}}, \ \text{where}  \\ 
& & \hspace{+3mm} z_{k_{0}}^{(0)} = z_{0} \in \widehat{\mathbb{C}}, \  \beta_{r} \in \left\{ 1, 2, \cdots, N \right\} \ \text{and}\ 1 \le k_{r} \le \delta_{\beta_{r}} \left( z_{k_{r - 1}}^{(-(r - 1))} \right)\ \text{for}\ r \in \mathbb{Z}_{+} \Big\}. 
\end{eqnarray*} 
Any point in $\mathscr{Q}^{\Gamma} \left( \widehat{\mathbb{C}} \right)$ is denoted by $\mathfrak{X}^{-} \left( z_{0}; \boldsymbol{\beta} \right)_{\boldsymbol{k}}$ where $\boldsymbol{\beta} = \left( \cdots, \beta_{2}, \beta_{1} \right)$ and $\boldsymbol{k} = \left( \cdots, k_{2}, k_{1} \right)$.
Similarly, we define the space of $n$ long and infinitely long forward orbits of points in $\widehat{\mathbb{C}}$ under the correspondence as follows. 
\begin{eqnarray*} 
\mathscr{P}_{n}^{\Gamma} \left( \widehat{\mathbb{C}} \right) & = & \Big\{ \left( z_{0}, z_{j_{1}}^{(1)}, \cdots, z_{j_{n}}^{(n)};\; \alpha_{1}, \cdots, \alpha_{n} \right)\ :\ \left(z_{j_{r-1}}^{(r-1)}, z_{j_{r}}^{(r)} \right) \in \Gamma_{\alpha_{r}} \ \text{where}\ z_{j_{0}}^{(0)} = z_{0} \in \widehat{\mathbb{C}}, \\ 
& & \hspace{+1cm} \alpha_{r} \in \left\{ 1, 2, \cdots, N \right\} \text{and}\ 1 \le j_{r} \le \lambda_{\alpha_{r}} \left( z_{j_{r - 1}}^{(r - 1)} \right)\ \text{for}\ 1 \le r \le n \Big\},\ \ \text{and} \\ 
\mathscr{P}^{\Gamma} \left( \widehat{\mathbb{C}} \right) & = & \Big\{ \left( z_{0}, z_{j_{1}}^{(1)}, z_{j_{2}}^{(2)}, \cdots;\; \alpha_{1}, \alpha_{2}, \cdots \right)\ :\ \left(z_{j_{r-1}}^{(r-1)}, z_{j_{r}}^{(r)} \right) \in \Gamma_{\alpha_{r}} \ \text{where}\ z_{j_{0}}^{(0)} = z_{0} \in \widehat{\mathbb{C}}, \\ 
& & \hspace{+1cm} \alpha_{r} \in \left\{ 1, 2, \cdots, N \right\}\ \text{and}\ 1 \le j_{r} \le \lambda_{\alpha_{r}} \left( z_{j_{r - 1}}^{(r - 1)} \right) \ \text{for}\ r \in \mathbb{Z}_{+} \Big\}.
\end{eqnarray*} 
The symbol $\lambda_{\alpha_{r}} \left( z_{j_{r - 1}}^{(r - 1)} \right)$ denotes the cardinality of the set $\digamma \left(z_{j_{r - 1}}^{(r - 1)} \right) \cap \pi_{2} (\Gamma_{\alpha_r})$. We denote any point in $\mathscr{P}_{n}^{\Gamma} \left( \widehat{\mathbb{C}} \right)$ by $\mathfrak{X}_{n}^{+} \left( z_{0}; \boldsymbol{\alpha} \right)_{\boldsymbol{j}}$ where $z_{0} \in \widehat{\mathbb{C}},\ \boldsymbol{\alpha} = \left( \alpha_{1}, \cdots, \alpha_{n} \right)$ and $\boldsymbol{j} = \left( j_{1}, \cdots, j_{n} \right)$ and any point in $\mathscr{P}^{\Gamma} \left( \widehat{\mathbb{C}} \right)$ by $\mathfrak{X}^{+} \left( z_{0}; \boldsymbol{\alpha} \right)_{\boldsymbol{j}}$ where $z_{0} \in \widehat{\mathbb{C}},\ \boldsymbol{\alpha} = \left( \alpha_{1}, \alpha_{2}, \cdots \right)$ and $\boldsymbol{j} = \left( j_{1}, j_{2}, \cdots \right)$. 

We now define a shift map $\sigma^{\Gamma} : \mathscr{P}^{\Gamma} \left( \widehat{\mathbb{C}} \right) \longrightarrow \mathscr{P}^{\Gamma} \left( \widehat{\mathbb{C}} \right)$ given by 
\[ \sigma^{\Gamma} \left( \mathfrak{X}^{+} \left( z_{0}; \boldsymbol{\alpha} \right)_{\boldsymbol{j}} \right)\ \ =\ \ \left( z_{j_{1}}^{(1)}, z_{j_{2}}^{(2)}, \cdots;\; \alpha_{2}, \cdots \right)\ \ =\ \ \mathfrak{X}^{+} \left( z_{j_{1}}^{(1)}; \sigma \boldsymbol{\alpha} \right)_{\sigma \boldsymbol{j}}, \] 
where by prefixing $\sigma$ to any infinite-lettered word $\boldsymbol{\alpha}$ or $\boldsymbol{j}$, we mean the infinite-lettered word that we obtain from $\boldsymbol{\alpha}$ or $\boldsymbol{j}$ as appropriate, by dropping its first letter. Also we use the notation $\sigma^{m} \boldsymbol{\alpha}$ and $\sigma^{m} \boldsymbol{j}$ to denote $\sigma \left( \sigma^{m - 1} \boldsymbol{\alpha} \right)$ and $\sigma \left( \sigma^{m - 1} \boldsymbol{j} \right)$ respectively, for any $m > 1$. 

We define projection maps $\Pi_{(r, n)}^{-} : \mathscr{Q}_{n}^{\Gamma} \left( \widehat{\mathbb{C}} \right) \longrightarrow \widehat{\mathbb{C}}$ given by $\Pi_{(r, n)}^{-}  \left( \mathfrak{X}^{-}_{n} \left( z_{0}; \boldsymbol{\beta} \right)_{\boldsymbol{k}} \right) = z_{k_{r}}^{(-r)}$ for $0 \le r \le n$ and $\Pi_{r}^{-} : \mathscr{Q}^{\Gamma} \left( \widehat{\mathbb{C}} \right) \longrightarrow \widehat{\mathbb{C}}$ given by $\Pi_{r}^{-}  \left( \mathfrak{X}^{-} \left( z_{0}; \boldsymbol{\beta} \right)_{\boldsymbol{j}} \right) = z_{k_{r}}^{(-r)}$ for $r \in \mathbb{Z}_{+} \cup \{ 0 \}$.

One can observe that $\mathscr{P}^{\Gamma} \left( \widehat{\mathbb{C}} \right)$ becomes a compact metric space when equipped with the metric 
\[ d \left( \mathfrak{X}^{+} \left( z_{0}; \boldsymbol{\alpha} \right)_{\boldsymbol{j}},\; \mathfrak{X}^{+} \left( w_{0}; \boldsymbol{\gamma} \right)_{\boldsymbol{l}} \right)\ =\ \max \left\{ \sum\limits_{n\, \ge\, 0} \frac{1}{2^{n}} \frac{\rho \left( z_{j_{n}}^{(n)},\, w_{l_{n}}^{(n)} \right)}{\left( 1 + \rho \left( z_{j_{n}}^{(n)},\, w_{l_{n}}^{(n)} \right) \right)}\ ,\ \frac{1}{2^{\min \left\{ n \ge 1 \; \mid\; \alpha_{n} \ne \gamma_{n} \right\}}} \right\}. \] 
This also results in $\sigma^{\Gamma}$ being a continuous map. 

Recall the set $\Omega_{{\rm DS}}$ which was introduced in Section \ref{sec:one}. As we mentioned before, we wish to study the restriction of $\digamma$ on $\Omega_{{\rm DS}}$. By an abuse of notation, we denote this restriction by $\digamma$ as well, for the sake of convenience. The study of this restriction is possible due to the following theorem due to Londhe, which we have paraphrased to suit our requirements.

\begin{theorem}\cite{lon:2022}
    Let $\digamma$ be a holomorphic correspondence on $\widehat{\mathbb{C}}$ and let $\Omega_{{\rm DS}}$ denote the support of the Dinh-Sibony measure of $\digamma$. Then we have $\digamma^{\dagger} (\Omega_{{\rm DS}}) \subset \Omega_{{\rm DS}}$. Further, for any $z \in \Omega_{{\rm DS}}$ we have $\digamma (z) \cap \Omega_{{\rm DS}} \neq \emptyset$.
\end{theorem} 

Thus, it makes sense to study the space of all infinitely long forward orbits under $\Gamma$ that lie entirely inside $\Omega_{{\rm DS}}$. On that note, we define $\mathscr{P}^{\Gamma}_{{\rm inv}} \left( \Omega_{{\rm DS}} \right) = \mathscr{P}^{\Gamma} \left( \widehat{\mathbb{C}} \right) \cap \left( \Omega_{{\rm DS}}^{\mathbb{Z}_{+}} \times \big\{ 1, 2, \cdots, N \big\}^{\mathbb{Z}_{+}} \right)$. This is a compact subset of $\mathscr{P}^{\Gamma} \left( \widehat{\mathbb{C}} \right)$. Further, the projection $\Pi_{0} : \mathscr{P}^{\Gamma}_{{\rm inv}} \left( \Omega_{{\rm DS}} \right) \longrightarrow \Omega_{{\rm DS}}$ is defined as $\Pi_{0} (\mathfrak{X}^{+} \left( z_{0}; \boldsymbol{\alpha} \right)_{\boldsymbol{j}}) = z_0$.

\section{Measure theoretic aspects of correspondences}
\label{sec:three}

In this section, we briefly set up the measure theoretic preliminaries required in our work. We also explore some of the latest results related to ergodic theory of holomorphic correspondences. The main sources of content in this section are \cite{ss:2025, ss:2026, bks:2026}. Studies on thermodynamic formalism in the dynamics of holomorphic correspondences date back to the works \cite{ds:2008, bs:2021}. The concept of pressure was first introduced in \cite{ss:2025} and that of measure theoretic entropy was first studied in \cite{ss:2026}. A theory of entropy and pressure for correspondences that is different from what we work with was developed in \cite{llz:2023}. They achieve various results on thermodynamic formalism for correspondences (set-valued maps), while we focus only on holomorphic correspondences. This idea lets us use the analyticity available in our setting, in order to prove results from a more complex dynamics perspective. From this section onwards, unless specified otherwise, we work with the correspondence $\digamma$ being restricted to $\Omega_{{\rm DS}}$. 

We denote by $\mathcal{M} (X)$ the space of all Borel probability measures defined on a compact metric space $X$. Note that the projection $\Pi_{0}$ mentioned in Section \ref{sec:two} induces a pushforward map $(\Pi_{0})_{*}: \mathcal{M} \left( \mathscr{P}^{\Gamma}_{{\rm inv}} \left( \Omega_{{\rm DS}} \right) \right) \longrightarrow \mathcal{M} (\Omega_{{\rm DS}})$. Define the space $\mathscr{S}^{\Gamma} := \left( \Pi_{0} \right)_{*} \left( \mathcal{M} \left( \mathscr{P}^{\Gamma}_{{\rm inv}} \left( \Omega_{{\rm DS}} \right), \sigma^{\Gamma} \right) \right)$, where $\mathcal{M} \left( \mathscr{P}^{\Gamma}_{{\rm inv}} \left( \Omega_{{\rm DS}} \right), \sigma^{\Gamma} \right)$ denotes the space of all $\sigma^{\Gamma}$- invariant measures in $\mathcal{M} \left( \mathscr{P}^{\Gamma}_{{\rm inv}} \left( \Omega_{{\rm DS}} \right) \right)$. For a measure $\nu \in \mathscr{S}^{\Gamma}$ with $\mu \in \mathcal{M} \left( \mathscr{P}^{\Gamma}_{{\rm inv}} \left( \Omega_{{\rm DS}} \right), \sigma^{\Gamma} \right)$ satisfying $(\Pi_{0})_{*} \mu = \nu$, the authors in \cite{ss:2026} defined the concept of intermediate entropy denoted by $h_{(\nu, \mu)} (\Gamma)$ and went on to define the entropy of $\Gamma$ with respect to $\nu$. The quantity $h_{(\nu, \mu)} (\Gamma)$ was defined using the classical partition construction using the projection maps and then by taking the supremum over all partitions. Readers may find all the necessary details regarding this in \cite{ss:2026}. 

\begin{definition}
\label{defn:entcorr}
    The entropy of $\Gamma$ with respect to a measure $\nu \in \mathscr{S}^{\Gamma}$ is defined as 
    \[ h_{\nu} (\Gamma)\ \ =\ \ \sup \left\{ h_{\left( \nu, \mu \right)} (\Gamma)\; :\; \mu \in  \mathcal{M} \left( \mathscr{P}^{\Gamma}_{{\rm inv}} \left( \Omega_{{\rm DS}} \right), \sigma^{\Gamma} \right)\ \text{such that}\ \nu = \left( \Pi_{0} \right)_{*} \mu \right\}. \]
\end{definition}

The notion of pressure of a potential $f \in \mathcal{C} (\Omega_{{\rm DS}}, \mathbb{R})$ under the dynamics of $\Gamma$, first studied in \cite{ss:2025}, is defined based on the classical idea of using separated and spanning sets. We shall denote the pressure of $f$ with respect to $\Gamma$ by ${\rm Pr} (\Gamma, f)$. We now state a few results from \cite{ss:2026} that are relevant to our work.

\begin{theorem}\cite{ss:2026}
\label{thm:entpresgamma}
    Let $\Gamma$ be a holomorphic correspondence on $\widehat{\mathbb{C}}$ and consider its restriction to the space $\Omega_{{\rm DS}}$. Let $\sigma^{\Gamma}$ denote the shift map on the space $\mathscr{P}^{\Gamma}_{{\rm inv}} \left( \Omega_{{\rm DS}} \right)$. Then we have,
    \begin{itemize}
        \item For any pair of measures $(\nu, \mu)$ such that $\nu \in \mathscr{S}^{\Gamma}$ with $\mu \in  \mathcal{M} \left( \mathscr{P}^{\Gamma}_{{\rm inv}} \left( \Omega_{{\rm DS}} \right), \sigma^{\Gamma} \right)$ satisfying $\nu = \left( \Pi_{0} \right)_{*} \mu$, we have $h_{(\nu, \mu)} (\Gamma) = h_{\mu} (\sigma^{\Gamma})$. Here $h_{\mu} (\sigma^{\Gamma})$ denotes the entropy of the measure $\mu$ with respect to the map $\sigma^{\Gamma}$. 
        \item For any $f \in \mathcal{C} (\Omega_{{\rm DS}}, \mathbb{R})$, we have ${\rm Pr} (\Gamma, f) = {\rm Pr} \left( \sigma^{\Gamma}, f \circ \Pi_{0} \right)$, where ${\rm Pr} \left( \sigma^{\Gamma}, f \circ \Pi_{0} \right)$ denotes the pressure of $f \circ \Pi_{0}$ with respect to $\sigma^{\Gamma}$.
    \end{itemize}
\end{theorem}

The next natural step after developing the concept of measure theoretic entropy and pressure for holomorphic correspondences was to prove a variational principle in the setting of holomorphic correspondences, which was achieved in \cite{ss:2026}. Their result is stated as follows.

\begin{theorem}\cite{ss:2026}
    Let $f \in \mathcal{C} (\Omega_{{\rm DS}}, \mathbb{R})$. Then we have,
    \begin{equation}
    \label{eqn:varprin}
    {\rm Pr} (\Gamma, f)\ \ =\ \ \sup\limits_{\nu\, \in\, \mathscr{S}^{\Gamma}} \left\{ h_{\nu} (\Gamma) + \int f\, \mathrm{d}\nu \right\}.
\end{equation}
\end{theorem}
In particular we have the topological entropy $h_{{\rm top}} (\Gamma) = \sup \left\{h_{\nu} (\Gamma) : \nu \in \mathscr{S}^{\Gamma} \right\}$. Any measure $\nu \in \mathscr{S}^{\Gamma}$ that achieves the supremum in Equation \eqref{eqn:varprin} is said to be an \emph{equilibrium state} for the function $f$. 

An interesting problem in thermodynamic formalism is to identify the interesting class of maps and potentials for which unique equilibrium state exists. While this theory is well explored in the case of single valued maps and in particular, in the study of rational maps with regards to complex dynamics, this question is relatively new in the study of holomorphic correspondences. The work \cite{bks:2026} addressed this question and identified an interesting class of correspondences to work with in this regard. In the study of single valued maps, typically one needs the map to be distance expanding, topologically transitive and open, while the potential should be H\"{o}lder continuous, in order for the existence of unique equilibrium state. Inspired by this, the authors in \cite{bks:2026} defined distance expanding correspondences and in fact showed that the conditions of openness and transitivity are automatically satisfied, when one restricts the correspondence to $\Omega_{{\rm DS}}$ (see Theorem 4.2 and Lemma 4.4 in \cite{bks:2026}), thereby generalising a similar study done in the case of a rational map restricted to its Julia set. Equipped with this, the authors proved results on the uniqueness of the equilibrium measure. We now state a definition and the main result from \cite{bks:2026}. All these results hold when the dynamics of the correspondence is restricted to $\Omega_{{\rm DS}}$.

\begin{definition}
Let $\digamma$ be a holomorphic correspondence on $\widehat{\mathbb{C}}$. The correspondence $\digamma$ is said to be \emph{distance expanding} on $\Omega_{{\rm DS}}$ if there exists $\lambda > 1,\ n \in \mathbb{Z}_{+}$ and $\eta > 0$ such that for all $z, w \in \Omega_{{\rm DS}}$ with $\rho (z, w) \le \eta$, we have, 
\[ \displaystyle{\inf \big\{ \rho (z_{1}, w_{1}) : z_{1} \in \digamma^{n} (z) \cap \Omega_{{\rm DS}}, w_{1} \in \digamma^{n} (w) \cap \Omega_{{\rm DS}} \big\} \ge \lambda\, \rho(z, w)}. \]
\end{definition}

\begin{theorem} \cite{bks:2026}
\label{thm:uniqueness} 
Let $\digamma$ be a holomorphic correspondence which is distance expanding on $\Omega_{{\rm DS}}$. Then, any H\"{o}lder continuous function $f \in \mathcal{C} \left( \Omega_{{\rm DS}}, \mathbb{R} \right)$ has a unique equilibrium state. 
\end{theorem}

\section{The bijection between $\mathscr{Q}_{n}^{\Gamma} \left( z_{0} \right)$ and $\left(\sigma^{\Gamma}\right)^{-n} \left(\mathfrak{X}^{+} \left( z_{0}; \boldsymbol{\alpha} \right)_{\boldsymbol{j}}\right)$}
\label{sec:four}

In this section, we state and prove some results that shall be useful in the proof of the main results, that shall appear in the next section. We fix an arbitrary point $\mathfrak{X}^{+} \left( z_{0}; \boldsymbol{\alpha} \right)_{\boldsymbol{j}} = \left( z_{0}, z_{j_{1}}^{(1)}, z_{j_{2}}^{(2)}, \cdots;\; \alpha_{1}, \alpha_{2}, \cdots \right) \in \mathscr{P}^{\Gamma}_{{\rm inv}} \left( \Omega_{{\rm DS}} \right)$ and work with it whenever necessary. We shall be interested in backward orbits for the most part, but fixing this forward path is needed in some of the calculations. However, even as we work with this fixed forward path, note that the results can be obtained by choosing any arbitrary fixed forward path and hence do not depend on the choice of $\mathfrak{X}^{+} \left( z_{0}; \boldsymbol{\alpha} \right)_{\boldsymbol{j}}$. 

An important observation that we make at this stage is that, for any given $z_0 \in \Omega_{{\rm DS}}$, the sets $\mathscr{Q}_{n}^{\Gamma} \left( z_{0} \right)$ and $\left(\sigma^{\Gamma}\right)^{-n} \left(\mathfrak{X}^{+} \left( z_{0}; \boldsymbol{\alpha} \right)_{\boldsymbol{j}}\right)$ are in a bijective correspondence. Indeed one simply maps $\mathfrak{X}^{+} \left( z_{k_{n}}^{(-n)}; \boldsymbol{\beta}\boldsymbol{\alpha} \right)_{\boldsymbol{k}\boldsymbol{j}} =  \left( z_{k_{n}}^{(-n)}, \cdots, z_{k_{1}}^{(-1)}, z_{0}, z_{j_{1}}^{(1)}, \cdots;\; \beta_{n}, \cdots, \beta_{1}, \alpha_{1}, \cdots \right) \in \left(\sigma^{\Gamma}\right)^{-n} \left(\mathfrak{X}^{+} \left( z_{0}; \boldsymbol{\alpha} \right)_{\boldsymbol{j}}\right)$ to the point $\mathfrak{X}_{n}^{-} \left( z_{0}; \boldsymbol{\beta} \right)_{\boldsymbol{k}} = \left( z_{k_{n}}^{(-n)}, \cdots, z_{k_{1}}^{(-1)}, z_{0};\; \beta_{n}, \cdots, \beta_{1} \right) \in \mathscr{Q}_{n}^{\Gamma} \left( z_{0} \right)$. 

As a result of this, we have the following consequences which shall be used in various proofs in the paper. For $\mathfrak{X}^{+} \left( z_{k_{n}}^{(-n)}; \boldsymbol{\beta}\boldsymbol{\alpha} \right)_{\boldsymbol{k}\boldsymbol{j}} \in \left(\sigma^{\Gamma}\right)^{-n} \left(\mathfrak{X}^{+} \left( z_{0}; \boldsymbol{\alpha} \right)_{\boldsymbol{j}}\right)$, if we set 
\[
\mu_{\mathfrak{X}^{+} \left( z_{k_{n}}^{(-n)}; \boldsymbol{\beta}\boldsymbol{\alpha} \right)_{\boldsymbol{k}\boldsymbol{j}}} \ \ = \ \  \displaystyle \frac{1}{n} \sum_{r \, =\, 0}^{n - 1} \delta_{\left(\sigma^{\Gamma}\right)^{r} \left(\mathfrak{X}^{+} \left( z_{k_{n}}^{(-n)}; \boldsymbol{\beta}\boldsymbol{\alpha} \right)_{\boldsymbol{k}\boldsymbol{j}}\right)} \ \ \in \ \ \mathcal{M} \left( \mathscr{P}^{\Gamma}_{{\rm inv}} \left( \Omega_{{\rm DS}} \right) \right),
\]
then a simple calculation shows that, $\left (\Pi_{0} \right)_{*} \left(\mu_{\mathfrak{X}^{+} \left( z_{k_{n}}^{(-n)}; \boldsymbol{\beta}\boldsymbol{\alpha} \right)_{\boldsymbol{k}\boldsymbol{j}}}\right) = \mu_{\mathfrak{X}_{n}^{-} \left( z_{0}; \boldsymbol{\beta} \right)_{\boldsymbol{k}}}$, where the quantity on the right hand side was defined in Equation \eqref{eqn:muzee}. Thus, for any $\mathcal{S} \subset \mathcal{M} (\Omega_{{\rm DS}})$, we have $\mu_{\mathfrak{X}_{n}^{-} \left( z_{0}; \boldsymbol{\beta} \right)_{\boldsymbol{k}}} \in \mathcal{S}$ if and only if $\mu_{\mathfrak{X}^{+} \left( z_{k_{n}}^{(-n)}; \boldsymbol{\beta}\boldsymbol{\alpha} \right)_{\boldsymbol{k}\boldsymbol{j}}} \in ((\Pi_{0})_{*})^{-1} (\mathcal{S})$.

For $f \in \mathcal{C} (\Omega_{{\rm DS}}, \mathbb{R})$, by denoting $F = f \circ \Pi_{0}$, we have,
\begin{equation}
\label{eqn:ergsumeq}
\mathcal{Z}_{n}^{-} (f) \left( \mathfrak{X}_{n}^{-} \left( z_{0}; \boldsymbol{\beta} \right)_{\boldsymbol{k}} \right) = \sum \limits_{r \, = \, 1}^{n} f\left(\Pi_{(r, n)}^{-} \left( \mathfrak{X}_{n}^{-} \left( z_{0}; \boldsymbol{\beta} \right)_{\boldsymbol{k}} \right)\right) = \sum \limits_{r \, = \, 0}^{n - 1} F \left(\left(\sigma^{\Gamma}\right)^{r} \left(\mathfrak{X}^{+} \left( z_{k_{n}}^{(-n)}; \boldsymbol{\beta}\boldsymbol{\alpha} \right)_{\boldsymbol{k}\boldsymbol{j}}\right)\right).
\end{equation}
Note that, the last expression in Equation \eqref{eqn:ergsumeq} is the ergodic sum of the potential $F$ with respect to $\sigma^{\Gamma}$, evaluated at $\mathfrak{X}^{+} \left( z_{k_{n}}^{(-n)}; \boldsymbol{\beta}\boldsymbol{\alpha} \right)_{\boldsymbol{k}\boldsymbol{j}}$, which is denoted by $F^{n} \left(\mathfrak{X}^{+} \left( z_{k_{n}}^{(-n)}; \boldsymbol{\beta}\boldsymbol{\alpha} \right)_{\boldsymbol{k}\boldsymbol{j}}\right)$. Further note that, if $f$ satisfies ${\rm Pr} (\Gamma, f) > \sup f$, then we also have that ${\rm Pr} (\sigma^{\Gamma}, F) > \sup F$. We are now ready to state a few results.

\begin{lemma}
\label{lem:limsuppr}
    Let $\Gamma$ be a holomorphic correspondence on $\widehat{\mathbb{C}}$ and let $f \in \mathcal{C} (\Omega_{{\rm DS}}, \mathbb{R})$. Fix $z_0 \in \Omega_{{\rm DS}}$. Then, we have $\displaystyle{\limsup_{n\, \to\, \infty} \frac{1}{n} \log \left( \sum_{\mathfrak{X}_{n}^{-} \left( z_{0}; \boldsymbol{\beta} \right)_{\boldsymbol{k}}\, \in\, \mathscr{Q}_{n}^{\Gamma} \left( z_{0} \right)} e^{\mathcal{Z}_{n}^{-} (f) \left( \mathfrak{X}_{n}^{-} \left( z_{0}; \boldsymbol{\beta} \right)_{\boldsymbol{k}} \right)} \right) \le {\rm Pr} \left( \Gamma, f \right)}.$  
\end{lemma}

\begin{proof}
    The proof of this lemma follows essentially from the ideas in \cite{ss:2025}. Note that, for any $z_0 \in \Omega_{{\rm DS}}$, one can find an $\epsilon_{z_0} > 0$ such that $\mathscr{Q}_{n}^{\Gamma} \left( z_{0} \right)$ is $(n, \epsilon_{z_{0}})$-separated, according to the definition in \cite{ss:2025}. This alongwith the definition of pressure associated to $f$ with respect to $\Gamma$ ensures that the result holds.
\end{proof}

Now we state a result which is a simple consequence of Proposition 3.4.3 in \cite{pu:2010}. Another fact that is used in proving the following lemma is that under the hypothesis that $\digamma$ is a distance expanding correspondence, the associated shift $\sigma^{\Gamma}$ defined on $\mathscr{P}^{\Gamma}_{{\rm inv}} \left( \Omega_{{\rm DS}} \right)$ is topologically transitive and also distance expanding (see Lemma 4.4 and Proposition 4.6 in \cite{bks:2026}).

\begin{lemma}
\label{lem:phinpres}
    Let $\digamma$ be a holomorphic correspondence which is distance expanding on $\Omega_{{\rm DS}}$. Let $F \in \mathcal{C} (\mathscr{P}^{\Gamma}_{{\rm inv}} \left( \Omega_{{\rm DS}} \right), \mathbb{R})$ be H\"{o}lder continuous. Then, for any $\mathfrak{X}^{+} \left( z_{0}; \boldsymbol{\alpha} \right)_{\boldsymbol{j}} \in \mathscr{P}^{\Gamma}_{{\rm inv}} \left( \Omega_{{\rm DS}} \right)$ we have,
    \[
    \lim \limits_{n \, \to \, \infty} \frac{1}{n} \log \left( \sum \limits_{\mathfrak{X}^{+} \left( z_{k_{n}}^{(-n)}; \boldsymbol{\beta}\boldsymbol{\alpha} \right)_{\boldsymbol{k}\boldsymbol{j}} \, \in \, \left(\sigma^{\Gamma}\right)^{-n} \left(\mathfrak{X}^{+} \left( z_{0}; \boldsymbol{\alpha} \right)_{\boldsymbol{j}}\right)}  e^{F^{n} \left(\mathfrak{X}^{+} \left( z_{k_{n}}^{(-n)}; \boldsymbol{\beta}\boldsymbol{\alpha} \right)_{\boldsymbol{k}\boldsymbol{j}}\right)} \right) = {\rm Pr} (\sigma^{\Gamma}, F).
    \]
\end{lemma}

Further, using the fact that $\left (\sigma^{\Gamma} \right)^{-n} (\mathfrak{X}^{+} \left( z_{0}; \boldsymbol{\alpha} \right)_{\boldsymbol{j}})$ is an $(n, \epsilon)$-separated set for a suitable $\epsilon$ (see \cite{pu:2010}) and the definition of the pressure of a continuous potential with respect to $\sigma^{\Gamma}$, we have for $F, G \in \mathcal{C} (\mathscr{P}^{\Gamma}_{{\rm inv}} \left( \Omega_{{\rm DS}} \right), \mathbb{R})$,
\begin{equation}
\label{eqn:phinpresplus}
    \limsup \limits_{n \, \to \, \infty} \frac{1}{n} \log \left( \sum \limits_{\mathfrak{X}^{+} \left( z_{k_{n}}^{(-n)}; \boldsymbol{\beta}\boldsymbol{\alpha} \right)_{\boldsymbol{k}\boldsymbol{j}} \, \in \, \left(\sigma^{\Gamma}\right)^{-n} \left(\mathfrak{X}^{+} \left( z_{0}; \boldsymbol{\alpha} \right)_{\boldsymbol{j}}\right)}  e^{(F+G)^{n} \left(\mathfrak{X}^{+} \left( z_{k_{n}}^{(-n)}; \boldsymbol{\beta}\boldsymbol{\alpha} \right)_{\boldsymbol{k}\boldsymbol{j}}\right)} \right) \le {\rm Pr} (\sigma^{\Gamma}, F + G).
\end{equation}

The following lemma is a result taken from \cite{wal:1982} (see Theorem 9.11 and Theorem 9.12 therein). This was also used in the proofs of results in \cite{ps:1996}. Recall that it was observed in \cite{bks:2026} that under the assumption that $\digamma$ is distance expanding, we have the entropy map $\mu \to h_{\mu} (\sigma^{\Gamma})$ to be upper semi-continuous on $\mathcal{M} \left( \mathscr{P}^{\Gamma}_{{\rm inv}} \left( \Omega_{{\rm DS}} \right), \sigma^{\Gamma} \right)$. 

\begin{lemma}
\label{lem:inflemma}
    Let $\digamma$ be a holomorphic correspondence which is distance expanding on $\Omega_{{\rm DS}}$. Consider the map $\sigma^{\Gamma}$ on $\mathscr{P}^{\Gamma}_{{\rm inv}} \left( \Omega_{{\rm DS}} \right)$. Then the following results hold.
    \begin{itemize}
        \item Suppose $\mu \in \mathcal{M} \left( \mathscr{P}^{\Gamma}_{{\rm inv}} \left( \Omega_{{\rm DS}} \right), \sigma^{\Gamma} \right)$, then 
        \[
        h_{\mu} (\sigma^{\Gamma}) = \inf \left\{ {\rm Pr} (\sigma^{\Gamma}, G) - \displaystyle \int G \, \mathrm{d}\mu : G \in \mathcal{C} (\mathscr{P}^{\Gamma}_{{\rm inv}} \left( \Omega_{{\rm DS}} \right), \mathbb{R})  \right\}.
        \]
        \item Suppose $\mu \in \mathcal{M} \left( \mathscr{P}^{\Gamma}_{{\rm inv}} \left( \Omega_{{\rm DS}} \right)\right) \setminus \mathcal{M} \left( \mathscr{P}^{\Gamma}_{{\rm inv}} \left( \Omega_{{\rm DS}} \right), \sigma^{\Gamma} \right)$, then
        \[
        \inf \left\{ {\rm Pr} (\sigma^{\Gamma}, G) - \displaystyle \int G \, \mathrm{d}\mu : G \in \mathcal{C} (\mathscr{P}^{\Gamma}_{{\rm inv}} \left( \Omega_{{\rm DS}} \right), \mathbb{R})  \right\} \le 0.
        \]
    \end{itemize}
\end{lemma}

We now define a rate function and prove some properties about it as that will play a key role in proving the large deviation result. We consider a holomorphic correspondence $\digamma$ which is distance expanding on $\Omega_{{\rm DS}}$. Let $F \in \mathcal{C} (\mathscr{P}^{\Gamma}_{{\rm inv}} \left( \Omega_{{\rm DS}} \right), \mathbb{R})$ be H\"{o}lder continuous. 

Define $Q_{F} : \mathcal{C} (\mathscr{P}^{\Gamma}_{{\rm inv}} \left( \Omega_{{\rm DS}} \right), \mathbb{R}) \longrightarrow \mathbb{R}$ given by $Q_{F} (G) = {\rm Pr} (\sigma^{\Gamma}, F + G) - {\rm Pr} (\sigma^{\Gamma}, F)$. For any measure $\mu \in \mathcal{M} \left( \mathscr{P}^{\Gamma}_{{\rm inv}} \left( \Omega_{{\rm DS}} \right) \right)$, define the rate function,
\[ I_{Q_{F}} (\mu)\ \ =\ \ \sup\limits_{G\, \in\, \mathcal{C} (\mathscr{P}^{\Gamma}_{{\rm inv}} \left( \Omega_{{\rm DS}} \right), \mathbb{R})} \left( \int G\, \mathrm{d}\mu - Q_{F} (G) \right). \] 

\begin{lemma} 
\label{lem:lsc} 
Consider a holomorphic correspondence $\digamma$ which is distance expanding on $\Omega_{{\rm DS}}$. Let $F : \mathscr{P}^{\Gamma}_{{\rm inv}} \left( \Omega_{{\rm DS}} \right) \longrightarrow \mathbb{R}$ be a H\"{o}lder continuous function whose unique equilibrium state is denoted by $\mu_{F} \in \mathcal{M} \left( \mathscr{P}^{\Gamma}_{{\rm inv}} \left( \Omega_{{\rm DS}} \right), \sigma^{\Gamma} \right)$. Also let ${\rm Pr} (\sigma^{\Gamma}, F) > \sup F$. Then, 
\begin{enumerate} 
\item $I_{Q_{F}} (\mu) > 0$ for all measures $\mu \ne \mu_{F}$ in $\mathcal{M} \left( \mathscr{P}^{\Gamma}_{{\rm inv}} \left( \Omega_{{\rm DS}} \right)\right)$, 
\item The map $\mu \longmapsto I_{Q_{F}} (\mu)$ is lower semi-continuous on $\mathcal{M} \left( \mathscr{P}^{\Gamma}_{{\rm inv}} \left( \Omega_{{\rm DS}} \right), \sigma^{\Gamma} \right)$ and it is bounded away from $0$ on $\mathcal{M} \left( \mathscr{P}^{\Gamma}_{{\rm inv}} \left( \Omega_{{\rm DS}} \right)\right) \setminus \mathcal{M} \left( \mathscr{P}^{\Gamma}_{{\rm inv}} \left( \Omega_{{\rm DS}} \right), \sigma^{\Gamma} \right)$. 
\end{enumerate} 
\end{lemma} 

\begin{proof} 
\begin{enumerate}
\item  For $\mu \in \mathcal{M} \left( \mathscr{P}^{\Gamma}_{{\rm inv}} \left( \Omega_{{\rm DS}} \right)\right)$,
\begin{eqnarray}
I_{Q_{F}} (\mu) & =  & \sup\limits_{G\, \in\, \mathcal{C} (\mathscr{P}^{\Gamma}_{{\rm inv}} \left( \Omega_{{\rm DS}} \right), \mathbb{R})} \left( \int G\, \mathrm{d}\mu - {\rm Pr} \left( \sigma^{\Gamma},\; F + G \right) + {\rm Pr} \left( \sigma^{\Gamma},\; F \right) \right) \nonumber \\
& = & \sup\limits_{G\, \in\, \mathcal{C} (\mathscr{P}^{\Gamma}_{{\rm inv}} \left( \Omega_{{\rm DS}} \right), \mathbb{R})} \left( \int (G - F)\, \mathrm{d}\mu - {\rm Pr} \left( \sigma^{\Gamma},\; G \right) + {\rm Pr} \left( \sigma^{\Gamma},\; F \right) \right) \nonumber \\
& = & \sup\limits_{G\, \in\, \mathcal{C} (\mathscr{P}^{\Gamma}_{{\rm inv}} \left( \Omega_{{\rm DS}} \right), \mathbb{R})} \left( \int G\, \mathrm{d}\mu - {\rm Pr} \left( \sigma^{\Gamma},\; G \right) \right) + {\rm Pr} \left( \sigma^{\Gamma},\; F \right) - \int F\, \mathrm{d}\mu  \nonumber \\
& = & - \inf\limits_{G\, \in\, \mathcal{C} (\mathscr{P}^{\Gamma}_{{\rm inv}} \left( \Omega_{{\rm DS}} \right), \mathbb{R})} \left({\rm Pr} \left( \sigma^{\Gamma},\; G \right) - \int G\, \mathrm{d}\mu \right) + {\rm Pr} \left( \sigma^{\Gamma},\; F \right) - \int F\, \mathrm{d}\mu  \nonumber
\end{eqnarray}

By Lemma \ref{lem:inflemma}, we have for $\mu \in \mathcal{M} \left( \mathscr{P}^{\Gamma}_{{\rm inv}} \left( \Omega_{{\rm DS}} \right), \sigma^{\Gamma} \right)$, 
\begin{equation}
\label{eqn:iqnuent}
   I_{Q_{F}} (\mu) = - h_{\mu} (\sigma^{\Gamma}) + {\rm Pr} \left( \sigma^{\Gamma},\; F \right) - \int F\, \mathrm{d}\mu.
\end{equation}

Now, in Equation \eqref{eqn:iqnuent} we use the fact that $\mu_F$ is the unique equilibrium state for $F$ to get $I_{Q_{F}} (\mu) > 0$ for all measures $\mu \in \mathcal{M} \left( \mathscr{P}^{\Gamma}_{{\rm inv}} \left( \Omega_{{\rm DS}} \right), \sigma^{\Gamma} \right)$ with $\mu \neq \mu_{F}$.

Again by Lemma \ref{lem:inflemma}, for $\mu \in \mathcal{M} \left( \mathscr{P}^{\Gamma}_{{\rm inv}} \left( \Omega_{{\rm DS}} \right)\right) \setminus \mathcal{M} \left( \mathscr{P}^{\Gamma}_{{\rm inv}} \left( \Omega_{{\rm DS}} \right), \sigma^{\Gamma} \right)$, 
\begin{eqnarray}
\label{eqn:iqnupos}
I_{Q_{F}} (\mu) & = & - \inf\limits_{G\, \in\, \mathcal{C} (\mathscr{P}^{\Gamma}_{{\rm inv}} \left( \Omega_{{\rm DS}} \right), \mathbb{R})} \left({\rm Pr} \left( \sigma^{\Gamma},\; G \right) - \int G\, \mathrm{d}\mu \right) + {\rm Pr} \left( \sigma^{\Gamma},\; F \right) - \int F\, \mathrm{d}\mu  \nonumber \\
& \ge & {\rm Pr} \left( \sigma^{\Gamma},\; F \right) - \int F\, \mathrm{d}\mu   \nonumber \\
I_{Q_{F}} (\mu) & \ge & {\rm Pr} \left( \sigma^{\Gamma},\; F \right) - \sup F > 0 
\end{eqnarray}
\item Since $\digamma$ is distance expanding, we have the entropy map $\mu \mapsto h_{\mu} (\sigma^{\Gamma})$ to be upper semi-continuous on $\mathcal{M} \left( \mathscr{P}^{\Gamma}_{{\rm inv}} \left( \Omega_{{\rm DS}} \right), \sigma^{\Gamma} \right)$. Thus, using the expression in Equation \eqref{eqn:iqnuent}, it follows that $\mu \mapsto I_{Q_{F}} (\mu)$ is lower semi-continuous on $\mathcal{M} \left( \mathscr{P}^{\Gamma}_{{\rm inv}} \left( \Omega_{{\rm DS}} \right), \sigma^{\Gamma} \right)$. Further, the expression in Equation \eqref{eqn:iqnupos} shows that $I_{Q_{F}}$ is bounded away from $0$ on $\mathcal{M} \left( \mathscr{P}^{\Gamma}_{{\rm inv}} \left( \Omega_{{\rm DS}} \right)\right) \setminus \mathcal{M} \left( \mathscr{P}^{\Gamma}_{{\rm inv}} \left( \Omega_{{\rm DS}} \right), \sigma^{\Gamma} \right)$. 
\end{enumerate} 
\end{proof}

\section{Proofs of main results}
\label{sec:five}

This section is dedicated to proving the main results of this paper, namely Theorem \ref{thm:main} and Theorem \ref{thm:muf}. Before we proceed to their proofs, we prove a lemma that shall be used in the proof of Theorem \ref{thm:main}. At this stage, we introduce a simplified notation that we use in the proof of the lemma, for the sake of brevity. 

For $H \in \mathcal{C} (\mathscr{P}^{\Gamma}_{{\rm inv}} \left( \Omega_{{\rm DS}} \right), \mathbb{R})$ and a compact subset $\mathcal{S}$ of $\mathcal{M} \left( \mathscr{P}^{\Gamma}_{{\rm inv}} \left( \Omega_{{\rm DS}} \right)\right)$, we denote 
\begin{eqnarray}
    \label{eqn:shortnotnepowernum}
    \Phi_{n, z_{0}}^{\mathcal{S}} (H) \ \ \ & = & \ \displaystyle \sum\limits_{\begin{minipage}{6cm} \begin{center} \scriptsize{$\mathfrak{X}^{+} \left( z_{k_{n}}^{(-n)}; \boldsymbol{\beta}\boldsymbol{\alpha} \right)_{\boldsymbol{k}\boldsymbol{j}} \, \in \, \left(\sigma^{\Gamma}\right)^{-n} \left(\mathfrak{X}^{+} \left( z_{0}; \boldsymbol{\alpha} \right)_{\boldsymbol{j}}\right)$} \\ \scriptsize{$\mu_{\mathfrak{X}^{+} \left( z_{k_{n}}^{(-n)}; \boldsymbol{\beta}\boldsymbol{\alpha} \right)_{\boldsymbol{k}\boldsymbol{j}}} \in \mathcal{S}$} \end{center} \end{minipage}} e^{H^{n} \left(\mathfrak{X}^{+} \left( z_{k_{n}}^{(-n)}; \boldsymbol{\beta}\boldsymbol{\alpha} \right)_{\boldsymbol{k}\boldsymbol{j}}\right)}  \\
    \label{eqn:eqn:shortnotnepowerdenom}
    \Phi_{n, z_{0}} (H) \ \ \ & = & \ \displaystyle \sum\limits_{\mathfrak{X}^{+} \left( z_{k_{n}}^{(-n)}; \boldsymbol{\beta}\boldsymbol{\alpha} \right)_{\boldsymbol{k}\boldsymbol{j}} \, \in \, \left(\sigma^{\Gamma}\right)^{-n} \left(\mathfrak{X}^{+} \left( z_{0}; \boldsymbol{\alpha} \right)_{\boldsymbol{j}}\right)} e^{H^{n} \left(\mathfrak{X}^{+} \left( z_{k_{n}}^{(-n)}; \boldsymbol{\beta}\boldsymbol{\alpha} \right)_{\boldsymbol{k}\boldsymbol{j}}\right)}
\end{eqnarray}

\begin{lemma} 
\label{lem:lf}
Let $\digamma$ be a holomorphic correspondence which is distance expanding on $\Omega_{{\rm DS}}$. Consider a H\"{o}lder continuous function $f: \Omega_{{\rm DS}} \longrightarrow \mathbb{R}$ and a point $z_0 \in \Omega_{{\rm DS}}$. Let $\mathcal{K} \subsetneq \mathcal{M} \left( \Omega_{{\rm DS}} \right)$ be a weak*-compact set. Then, there exists $I_{\mathcal{K}} > 0$ such that,
\[ \limsup_{n\, \to\, \infty} \frac{1}{n} \log \left( \dfrac{\displaystyle \sum\limits_{\begin{minipage}{3cm} \begin{center} \scriptsize{$\mathfrak{X}_{n}^{-} \left( z_{0}; \boldsymbol{\beta} \right)_{\boldsymbol{k}}\, \in\, \mathscr{Q}_{n}^{\Gamma} \left( z_{0} \right)$} \\ \scriptsize{$\mu_{\mathfrak{X}_{n}^{-} \left( z_{0}; \boldsymbol{\beta} \right)_{\boldsymbol{k}}}\, \in\, \mathcal{K}$} \end{center} \end{minipage}} e^{\mathcal{Z}_{n}^{-} (f) \left( \mathfrak{X}_{n}^{-} \left( z_{0}; \boldsymbol{\beta} \right)_{\boldsymbol{k}} \right)}}{\displaystyle \sum\limits_{\mathfrak{X}_{n}^{-} \left( z_{0}; \boldsymbol{\beta} \right)_{\boldsymbol{k}}\, \in\, \mathscr{Q}_{n}^{\Gamma} \left( z_{0} \right)} e^{\mathcal{Z}_{n}^{-} (f) \left( \mathfrak{X}_{n}^{-} \left( z_{0}; \boldsymbol{\beta} \right)_{\boldsymbol{k}} \right)}} \right)\ \ \le\ \ - I_{\mathcal{K}}. \]
\end{lemma} 

\begin{proof}
Consider the weak* compact subset $((\Pi_{0})_{*})^{-1} (\mathcal{K}) \subset \mathcal{M} \left( \mathscr{P}^{\Gamma}_{{\rm inv}} \left( \Omega_{{\rm DS}} \right)\right)$. We denote by $F$ the real valued H\"{o}lder continuous function $f \circ \Pi_{0}$ defined on $\mathscr{P}^{\Gamma}_{{\rm inv}} \left( \Omega_{{\rm DS}} \right)$. 
    
Set $I_{\mathcal{K}} = \inf \limits_{\mu \, \in \, ((\Pi_{0})_{*})^{-1} (\mathcal{K})} I_{Q_{F}} (\mu)$. Now let $\epsilon > 0$ and $\mu \in ((\Pi_{0})_{*})^{-1} (\mathcal{K})$. Then, there exists $G \in \mathcal{C} (\mathscr{P}^{\Gamma}_{{\rm inv}} \left( \Omega_{{\rm DS}} \right), \mathbb{R})$ such that $\displaystyle \int G \, \mathrm{d}\mu - Q_{F} (G) > I_{Q_{F}} (\mu) - \epsilon \ge I_{\mathcal{K}} - \epsilon$. Thus, we have
\[
((\Pi_{0})_{*})^{-1} (\mathcal{K}) \ \subset \bigcup \limits_{G \, \in \, \mathcal{C} \left(\mathscr{P}^{\Gamma}_{{\rm inv}} \left( \Omega_{{\rm DS}} \right), \mathbb{R}\right)} \left\{ \mu \in \mathcal{M} \left( \mathscr{P}^{\Gamma}_{{\rm inv}} \left( \Omega_{{\rm DS}} \right)\right) : \int G \, \mathrm{d} \mu - Q_{F} (G) > I_{\mathcal{K}} - \epsilon \right\}.
\]
Since $((\Pi_{0})_{*})^{-1} (\mathcal{K})$ is weak*-compact, there exists $G_j \in \mathcal{C} (\mathscr{P}^{\Gamma}_{{\rm inv}} \left( \Omega_{{\rm DS}} \right), \mathbb{R})$ with $1 \le j \le k$ such that 
\[
((\Pi_{0})_{*})^{-1} (\mathcal{K}) \subset \bigcup \limits_{j = 1}^{k} \left\{ \mu \in \mathcal{M} \left( \mathscr{P}^{\Gamma}_{{\rm inv}} \left( \Omega_{{\rm DS}} \right)\right) : \int G_j \, \mathrm{d} \mu - Q_{F} (G_j) > I_{\mathcal{K}} - \epsilon \right\}.
\]
Recall the arbitrary fixed forward orbit $\mathfrak{X}^{+} \left( z_{0}; \boldsymbol{\alpha} \right)_{\boldsymbol{j}} \in \mathscr{P}^{\Gamma}_{{\rm inv}} \left( \Omega_{{\rm DS}} \right)$ and the discussion at the beginning of Section \ref{sec:four}. Now suppose $\mathfrak{X}_{n}^{-} \left( z_{0}; \boldsymbol{\beta} \right)_{\boldsymbol{k}} \in \mathscr{Q}_{n}^{\Gamma} \left( z_{0} \right)$ is such that $\mu_{\mathfrak{X}_{n}^{-} \left( z_{0}; \boldsymbol{\beta} \right)_{\boldsymbol{k}}} \in \mathcal{K}$, then we have the corresponding point $\mathfrak{X}^{+} \left( z_{k_{n}}^{(-n)}; \boldsymbol{\beta}\boldsymbol{\alpha} \right)_{\boldsymbol{k}\boldsymbol{j}} \in \left(\sigma^{\Gamma}\right)^{-n} \left(\mathfrak{X}^{+} \left( z_{0}; \boldsymbol{\alpha} \right)_{\boldsymbol{j}}\right)$ to satisfy $\mu_{\mathfrak{X}^{+} \left( z_{k_{n}}^{(-n)}; \boldsymbol{\beta}\boldsymbol{\alpha} \right)_{\boldsymbol{k}\boldsymbol{j}}} \in ((\Pi_{0})_{*})^{-1} (\mathcal{K})$. Thus, there exists $1 \le j \le k$ and $G_{j} \in \mathcal{C} \left(\mathscr{P}^{\Gamma}_{{\rm inv}} \left( \Omega_{{\rm DS}} \right), \mathbb{R}\right)$ satisfying $\displaystyle \int G_{j} \, \mathrm{d} \mu_{\mathfrak{X}^{+} \left( z_{k_{n}}^{(-n)}; \boldsymbol{\beta}\boldsymbol{\alpha} \right)_{\boldsymbol{k}\boldsymbol{j}}} - \ Q_{F} (G_{j}) > I_{\mathcal{K}} - \epsilon$.

Thus, we get
\[
\int G_j \, \mathrm{d} \mu_{\mathfrak{X}^{+} \left( z_{k_{n}}^{(-n)}; \boldsymbol{\beta}\boldsymbol{\alpha} \right)_{\boldsymbol{k}\boldsymbol{j}}} - Q_{F} (G_j) = \frac{G_{j}^{n} \left(\mathfrak{X}^{+} \left( z_{k_{n}}^{(-n)}; \boldsymbol{\beta}\boldsymbol{\alpha} \right)_{\boldsymbol{k}\boldsymbol{j}}\right)}{n} - Q_{F} (G_j) >  I_{\mathcal{K}} - \epsilon,
\]

which gives $G_{j}^{n} \left(\mathfrak{X}^{+} \left( z_{k_{n}}^{(-n)}; \boldsymbol{\beta}\boldsymbol{\alpha} \right)_{\boldsymbol{k}\boldsymbol{j}}\right) - n (I_{\mathcal{K}} - \epsilon) -n Q_{F}(G_j) > 0$. 

Now observe that, 
\[
\limsup_{n\, \to\, \infty} \frac{1}{n} \log \left( \dfrac{\displaystyle \sum\limits_{\begin{minipage}{3cm} \begin{center} \scriptsize{$\mathfrak{X}_{n}^{-} \left( z_{0}; \boldsymbol{\beta} \right)_{\boldsymbol{k}}\, \in\, \mathscr{Q}_{n}^{\Gamma} \left( z_{0} \right)$} \\ \scriptsize{$\mu_{\mathfrak{X}_{n}^{-} \left( z_{0}; \boldsymbol{\beta} \right)_{\boldsymbol{k}}}\, \in\, \mathcal{K}$} \end{center} \end{minipage}} e^{\mathcal{Z}_{n}^{-} (f) \left( \mathfrak{X}_{n}^{-} \left( z_{0}; \boldsymbol{\beta} \right)_{\boldsymbol{k}} \right)}}{\displaystyle \sum\limits_{\mathfrak{X}_{n}^{-} \left( z_{0}; \boldsymbol{\beta} \right)_{\boldsymbol{k}}\, \in\, \mathscr{Q}_{n}^{\Gamma} \left( z_{0} \right)} e^{\mathcal{Z}_{n}^{-} (f) \left( \mathfrak{X}_{n}^{-} \left( z_{0}; \boldsymbol{\beta} \right)_{\boldsymbol{k}} \right)}} \right) = \limsup \limits_{n \, \to \, \infty} \dfrac{1}{n} \log \left(\dfrac{\Phi_{n, z_0}^{((\Pi_{0})_{*})^{-1} (\mathcal{K})} (F)}{\Phi_{n, z_0} (F)}\right)
\]
and hence it is enough to estimate the quantity on the right hand side of the above equation. From above discussion, we have,
\begin{eqnarray*}
    \Phi_{n, z_0}^{((\Pi_{0})_{*})^{-1} (\mathcal{K})} (F) & \le & \sum_{j\, =\, 1}^{k} \sum\limits_{\begin{minipage}{7cm} \begin{center} \scriptsize{$\mathfrak{X}^{+} \left( z_{k_{n}}^{(-n)}; \boldsymbol{\beta}\boldsymbol{\alpha} \right)_{\boldsymbol{k}\boldsymbol{j}} \, \in \, \left(\sigma^{\Gamma}\right)^{-n} \left(\mathfrak{X}^{+} \left( z_{0}; \boldsymbol{\alpha} \right)_{\boldsymbol{j}}\right)$} \\ \scriptsize{$\frac{G_{j}^{n} \left(\mathfrak{X}^{+} \left( z_{k_{n}}^{(-n)}; \boldsymbol{\beta}\boldsymbol{\alpha} \right)_{\boldsymbol{k}\boldsymbol{j}}\right)}{n} -  Q_{F} (G_j) >  (I_{\mathcal{K}} - \epsilon)$} \end{center} \end{minipage}} e^{F^{n} \left(\mathfrak{X}^{+} \left( z_{k_{n}}^{(-n)}; \boldsymbol{\beta}\boldsymbol{\alpha} \right)_{\boldsymbol{k}\boldsymbol{j}}\right)} \\  
    & \le & \sum \limits_{j\, =\, 1}^{k}\ e^{-n \left( Q_{F}(G_{j}) + I_{\mathcal{K}} - \epsilon \right)} \ \  \Phi_{n, z_{0}} (F + G_{j}) \\
    & \le & k \sup \limits_{1\, \le\, j\, \le\, k} \left \{  e^{-n \left( Q_{F}(G_{j}) + I_{\mathcal{K}} - \epsilon \right)} \ \  \Phi_{n, z_{0}} (F + G_{j}) \right \}.
\end{eqnarray*}

Thus we get,
\begin{eqnarray*}
    \frac{1}{n} \log \left(\dfrac{\Phi_{n, z_0}^{((\Pi_{0})_{*})^{-1} (\mathcal{K})} (F)}{\Phi_{n, z_0} (F)}\right) & \le & \frac{1}{n} \log \left (\frac{ k \sup \limits_{1 \le j \le k} \left \{  e^{-n \left( Q_{F}(G_{j}) + I_{\mathcal{K}} - \epsilon \right)} \ \  \Phi_{n, z_{0}} (F + G_{j}) \right \}} {\Phi_{n, z_{0}} (F)} \right) \\
    & \le & \frac{1}{n} \log k - \frac{1}{n} \log \Phi_{n, z_0} (F) \\
    & & \hspace{7mm} + \sup \limits_{1 \le j \le k} \left \{ -Q_{F}(G_{j}) - I_{\mathcal{K}} + \epsilon + \frac{1}{n} \log \Phi_{n, z_0} (F + G_j) \right \},
\end{eqnarray*}
and hence,
\begin{eqnarray*}
    \limsup_{n \, \to \, \infty} \frac{1}{n} \log \left(\dfrac{\Phi_{n, z_0}^{((\Pi_{0})_{*})^{-1} (\mathcal{K})} (F)}{\Phi_{n, z_0} (F)}\right) & \le & \limsup \limits_{n \, \to \, \infty} \left ( \sup \limits_{1 \le j \le k} \Bigg \{ -Q_{F}(G_{j}) - I_{\mathcal{K}} + \epsilon   \right. \\
    & & \hspace{7mm}  \left. + \frac{1}{n} \log \Phi_{n, z_0} (F + G_j) \Bigg \} \right ) - \liminf \limits_{n \, \to \, \infty} \frac{1}{n} \log \Phi_{n, z_0} (F) \\
    & \le & \sup \limits_{1 \le j \le k} \left \{ -Q_{F}(G_{j}) - I_{\mathcal{K}} + \epsilon + \limsup \limits_{n \, \to \, \infty} \frac{1}{n} \log \Phi_{n, z_0} (F + G_{j})\right\} \\
    & & \hspace{10mm} - \liminf \limits_{n \, \to \, \infty} \frac{1}{n} \log \Phi_{n, z_0} (F).
\end{eqnarray*}
At this stage, an application of Lemma \ref{lem:phinpres} and Equation \eqref{eqn:phinpresplus} gives, 
\begin{eqnarray*}
    \limsup_{n \, \to \, \infty} \frac{1}{n} \log \left(\dfrac{\Phi_{n, z_0}^{((\Pi_{0})_{*})^{-1} (\mathcal{K})} (F)}{\Phi_{n, z_0} (F)}\right) & \le &  \sup \limits_{1 \le j \le k} \left \{ -Q_{F}(G_{j}) - I_{\mathcal{K}} + \epsilon + {\rm Pr} (\sigma^{\Gamma}, F + G_{j}) - {\rm Pr} (\sigma^{\Gamma}, F) \right\} \\
    & \le & - I_{\mathcal{K}} + \epsilon.
\end{eqnarray*}
As $\epsilon > 0$ is arbitrary, the result follows by taking $\epsilon \, \to \, 0$. 
\end{proof}

Now we present the proofs of the main results of the paper.

\begin{proof} [of Theorem \ref{thm:main}]
    Let $\mathcal{V} \subsetneq \mathcal{M} \left( \Omega_{{\rm DS}} \right)$ be an open neighbourhood of $\nu_{f}$, the unique equilibrium state of the potential $f$ and let $\nu_f = (\Pi_{0})_{*} \mu_{F}$. Set $\mathcal{K} = \mathcal{M} \left( \Omega_{{\rm DS}} \right) \setminus \mathcal{V}$. Note that as $\nu_{f} \notin \mathcal{K}$, we have $\mu_{F} \notin ((\Pi_{0})_{*})^{-1} (\mathcal{K})$. From Lemma \ref{lem:lsc}, it follows that $I_{\mathcal{K}} > 0$. Now, applying Lemma \ref{lem:lf} for the set $\mathcal{K} = \mathcal{M} \left({\Omega_{{\rm DS}}}\right) \setminus \mathcal{V}$ and then simplifying the inequality we get as a result, completes the proof of this theorem.
\end{proof}

For the sake of brevity, we adopt the notation given below in the proof that follows. Let $\mathcal{S} \subset \mathcal{M} \left( \Omega_{{\rm DS}} \right)$, $z_{0}$ be a point in $\Omega_{{\rm DS}}$ and $f, g \in \mathcal{C} (\Omega_{{\rm DS}}, \mathbb{R})$. 

\begin{eqnarray}
     \Psi_{n, z_{0}}^{\mathcal{S}} (f, g) \ \ \ & = & \sum\limits_{\begin{minipage}{3cm} \begin{center} \scriptsize{$ \mathfrak{X}_{n}^{-} \left( z_{0}; \boldsymbol{\beta} \right)_{\boldsymbol{k}} \, \in\,  \mathscr{Q}_{n}^{\Gamma} \left( z_{0} \right)$} \\ \scriptsize{$\mu_{\mathfrak{X}_{n}^{-} \left( z_{0}; \boldsymbol{\beta} \right)_{\boldsymbol{k}}} \, \in \,  \mathcal{S}$} \end{center} \end{minipage}} e^{\mathcal{Z}_{n}^{-} (f) \left( \mathfrak{X}_{n}^{-} \left( z_{0}; \boldsymbol{\beta} \right)_{\boldsymbol{k}} \right)} \, \frac{\mathcal{Z}_{n}^{-} (g) \left( \mathfrak{X}_{n}^{-} \left( z_{0}; \boldsymbol{\beta} \right)_{\boldsymbol{k}} \right)}{n} \\
     \Psi_{n, z_{0}} (f) & = & \sum\limits_{\mathfrak{X}_{n}^{-} \left( z_{0}; \boldsymbol{\beta} \right)_{\boldsymbol{k}}\, \in\, \mathscr{Q}_{n}^{\Gamma} \left( z_{0} \right)} e^{\mathcal{Z}_{n}^{-} (f) \left( \mathfrak{X}_{n}^{-} \left( z_{0}; \boldsymbol{\beta} \right)_{\boldsymbol{k}} \right)}
\end{eqnarray}

\begin{proof}[of Theorem \ref{thm:muf}] 
Let $\nu_n = \left( \dfrac{\sum\limits_{\mathfrak{X}_{n}^{-} \left( z_{0}; \boldsymbol{\beta} \right)_{\boldsymbol{k}}\, \in\, \mathscr{Q}_{n}^{\Gamma} \left( z_{0} \right)} e^{\mathcal{Z}_{n}^{-} (f) \left( \mathfrak{X}_{n}^{-} \left( z_{0}; \boldsymbol{\beta} \right)_{\boldsymbol{k}} \right)} \mu_{\mathfrak{X}_{n}^{-} \left( z_{0}; \boldsymbol{\beta} \right)_{\boldsymbol{k}}}}{\sum\limits_{\mathfrak{X}_{n}^{-} \left( z_{0}; \boldsymbol{\beta} \right)_{\boldsymbol{k}}\, \in\, \mathscr{Q}_{n}^{\Gamma} \left( z_{0} \right)} e^{\mathcal{Z}_{n}^{-} (f) \left( \mathfrak{X}_{n}^{-} \left( z_{0}; \boldsymbol{\beta} \right)_{\boldsymbol{k}} \right)}}\right)$. 
    
We prove that for $g \in \mathcal{C} (\Omega_{{\rm DS}}, \mathbb{R})$,  $\displaystyle{\int g \, \mathrm{d}\nu_n \to \int g \, \mathrm{d}\nu_{f}}$. To show this, for $\epsilon > 0$ we define the neighbourhood $\mathcal{V}$ of $\nu_{f}$ given by 
\[\displaystyle{\mathcal{V} := \left\{ \nu \in \mathcal{M} \left(\Omega_{{}\rm DS} \right) : \left| \int g \, \mathrm{d} \nu  - \int g \, \mathrm{d} \nu_f\right| < \epsilon, \ \forall g \, \in \, \mathcal{C} (\Omega_{{\rm DS}}, \mathbb{R})  \right\}}.\]

We have for $g \, \in \, \mathcal{C} (\Omega_{{\rm DS}}, \mathbb{R})$,
\begin{eqnarray}
\label{eqn:gnunone}
\int g \, \mathrm{d}\nu_{n} & = & \frac{1}{\Psi_{n, z_{0}} (f)} \bigintss \sum\limits_{\mathfrak{X}_{n}^{-} \left( z_{0}; \boldsymbol{\beta} \right)_{\boldsymbol{k}}\, \in\, \mathscr{Q}_{n}^{\Gamma} \left( z_{0} \right)} e^{\mathcal{Z}_{n}^{-} (f) \left( \mathfrak{X}_{n}^{-} \left( z_{0}; \boldsymbol{\beta} \right)_{\boldsymbol{k}} \right)} \, g \, \mathrm{d} \hspace{-1 mm} \left( \frac{1}{n} \sum_{r\, =\, 1}^{n} \delta_{\Pi_{(r,n)}^{-} \left( \mathfrak{X}_{n}^{-} \left( z_{0}; \boldsymbol{\beta} \right)_{\boldsymbol{k}} \right)}\right) \nonumber \\
& = & \frac{1}{\Psi_{n, z_{0}} (f)} \sum\limits_{\mathfrak{X}_{n}^{-} \left( z_{0}; \boldsymbol{\beta} \right)_{\boldsymbol{k}}\, \in\, \mathscr{Q}_{n}^{\Gamma} \left( z_{0} \right)}  e^{\mathcal{Z}_{n}^{-} (f) \left( \mathfrak{X}_{n}^{-} \left( z_{0}; \boldsymbol{\beta} \right)_{\boldsymbol{k}} \right)} \, \frac{\mathcal{Z}_{n}^{-} (g) \left( \mathfrak{X}_{n}^{-} \left( z_{0}; \boldsymbol{\beta} \right)_{\boldsymbol{k}} \right)}{n} \nonumber \\
& = & \frac{1}{\Psi_{n, z_{0}} (f)} \left( \Psi_{n, z_{0}}^{\mathcal{V}} (f, g) + \Psi_{n, z_{0}}^{\mathcal{V}^{c}} (f, g)  \right).
\end{eqnarray}
   
Note that the quantity $\dfrac{\mathcal{Z}_{n}^{-} (g) \left( \mathfrak{X}_{n}^{-} \left( z_{0}; \boldsymbol{\beta} \right)_{\boldsymbol{k}} \right)}{n}$ is bounded for any $g \in \mathcal{C} (\Omega_{{\rm DS}}, \mathbb{R})$. This and an application of Theorem \ref{thm:main} gives $\dfrac{1}{\Psi_{n, z_{0}} (f)} \left( \Psi_{n, z_{0}}^{\mathcal{V}^{c}} (f, g)  \right)  \le Cr^n$, for some $0 < r < 1$ and $C > 0$. 
    
In the case where $\mu_{\mathfrak{X}_{n}^{-} \left( z_{0}; \boldsymbol{\beta} \right)_{\boldsymbol{k}}} \in \mathcal{V}$, we have the inequality $\displaystyle{\left | \frac{\mathcal{Z}_{n}^{-} (g) \left( \mathfrak{X}_{n}^{-} \left( z_{0}; \boldsymbol{\beta} \right)_{\boldsymbol{k}} \right)}{n} - \int g \, \mathrm{d}\nu_f \right| < \epsilon}$, for any $g \in \mathcal{C} (\Omega_{{\rm DS}}, \mathbb{R})$. Thus, we get
\[
\Psi_{n, z_{0}}^{\mathcal{V}} (f, g) = \displaystyle \sum\limits_{\begin{minipage}{2cm} \begin{center} \scriptsize{\begin{eqnarray*} \mathfrak{X}_{n}^{-} \left( z_{0}; \boldsymbol{\beta} \right)_{\boldsymbol{k}} & \hspace{-5pt} \in \hspace{-5pt} & \mathscr{Q}_{n}^{\Gamma} \left( z_{0} \right) \\ \mu_{\mathfrak{X}_{n}^{-} \left( z_{0}; \boldsymbol{\beta} \right)_{\boldsymbol{k}}} & \hspace{-5pt} \in \hspace{-5pt} & \mathcal{V} \end{eqnarray*}} \end{center} \end{minipage}}  e^{\mathcal{Z}_{n}^{-} (f) \left( \mathfrak{X}_{n}^{-} \left( z_{0}; \boldsymbol{\beta} \right)_{\boldsymbol{k}} \right)} \left ( \int g \, \mathrm{d} \nu_f + E_{n} \left( \mathfrak{X}_{n}^{-} \left( z_{0}; \boldsymbol{\beta} \right)_{\boldsymbol{k}} \right) \right), \nonumber
\]
with the error term $E_{n} \left( \mathfrak{X}_{n}^{-} \left( z_{0}; \boldsymbol{\beta} \right)_{\boldsymbol{k}} \right)$ satisfying $\left |E_{n} \left( \mathfrak{X}_{n}^{-} \left( z_{0}; \boldsymbol{\beta} \right)_{\boldsymbol{k}} \right) \right| < \epsilon$.

Thus, finally Equation \eqref{eqn:gnunone} becomes,
\[
\int g \, \mathrm{d} \nu_n = O (r^n) + \frac{1}{\Psi_{n, z_{0}} (f)} \displaystyle \sum\limits_{\begin{minipage}{2cm} \begin{center} \scriptsize{\begin{eqnarray*} \mathfrak{X}_{n}^{-} \left( z_{0}; \boldsymbol{\beta} \right)_{\boldsymbol{k}} & \hspace{-5pt} \in \hspace{-5pt} & \mathscr{Q}_{n}^{\Gamma} \left( z_{0} \right) \\ \mu_{\mathfrak{X}_{n}^{-} \left( z_{0}; \boldsymbol{\beta} \right)_{\boldsymbol{k}}} & \hspace{-5pt} \in \hspace{-5pt} & \mathcal{V} \end{eqnarray*}} \end{center} \end{minipage}}  e^{\mathcal{Z}_{n}^{-} (f) \left( \mathfrak{X}_{n}^{-} \left( z_{0}; \boldsymbol{\beta} \right)_{\boldsymbol{k}} \right)} \left ( \int g \, \mathrm{d} \nu_f + E_{n} \left( \mathfrak{X}_{n}^{-} \left( z_{0}; \boldsymbol{\beta} \right)_{\boldsymbol{k}} \right) \right) \nonumber
\]
and hence, we get
   \[
    \limsup \limits_{n \to \infty} \int g \, \mathrm{d}\nu_n \le \int g \, \mathrm{d}\nu_f + \epsilon \hspace{3mm} \text{and} \hspace{3mm} \liminf \limits_{n \to \infty} \int g \, \mathrm{d}\nu_n \ge \int g \, \mathrm{d}\nu_f - \epsilon 
   \]
Letting $\epsilon \to 0$ completes the proof of the theorem.     
\end{proof}

\section{The Dinh-Sibony measure and ergodicity}
\label{sec:six}

This section is devoted to exploring more about the Dinh-Sibony measure associated to a holomorphic correspondence and how it is connected to Theorem \ref{thm:muf}. In this regard, we also give some comments connecting various other works in the theory of holomorphic correspondences. Firstly, we prove Corollary \ref{cor:dsisuniqueent}. 

\begin{proof}[of Corollary \ref{cor:dsisuniqueent}]
    Note that the unique measure of maximal entropy is the unique equilibrium state corresponding to the function $f \equiv 0$. Thus, to prove this corollary, it is enough to prove that, the measure $\nu_{0}$ which is obtained as the limit of the sequence $\nu_n = \dfrac{1}{d_{{\rm top}}^{n}} \displaystyle \sum\limits_{\mathfrak{X}_{n}^{-} \left( z_{0}; \boldsymbol{\beta} \right)_{\boldsymbol{k}}\, \in\, \mathscr{Q}_{n}^{\Gamma} \left( z_{0} \right)}  \mu_{\mathfrak{X}_{n}^{-} \left( z_{0}; \boldsymbol{\beta} \right)_{\boldsymbol{k}}}$ in Theorem \ref{thm:muf} is indeed the Dinh-Sibony measure. A simple calculation shows that, $\nu_n = \dfrac{1}{n} \displaystyle \sum \limits_{r \, = \, 1}^{n} \dfrac{1}{d_{{\rm top}}^{r}} (F^{r})^{*}\delta_{z_{0}}$, which converges to $\mu_{\digamma}$. Hence the result follows.
\end{proof}

This also explains how Theorem \ref{thm:muf} generalises the convergence in Equation \eqref{eqn:dscgt}. We now present an alternative proof of Corollary \ref{cor:dsisuniqueent} and as a consequence, mention a few comments related to the recent work \cite{tmss:2026}.

\begin{proof}[alternative proof of Corollary \ref{cor:dsisuniqueent}]
    As done in Section \ref{sec:four}, fix  $\mathfrak{X}^{+} \left( z_{0}; \boldsymbol{\alpha} \right)_{\boldsymbol{j}} \in \mathscr{P}^{\Gamma}_{{\rm inv}} \left( \Omega_{{\rm DS}} \right)$. From Proposition 4.6 in \cite{bks:2026}, we have that $\sigma^{\Gamma}$ is distance expanding on $\mathscr{P}^{\Gamma}_{{\rm inv}} \left( \Omega_{{\rm DS}} \right)$. Now, using Theorem 7.1 in \cite{lop:1994}, there exists a $\sigma^{\Gamma}$ invariant ergodic measure $\mu$ such that, 
    \begin{equation}
    \label{eqn:sigmadsseq}
        \mu_{n} = \frac{1}{d_{{\rm top}}^{n}} \sum \limits_{\mathfrak{X}^{+} \left( z_{k_{n}}^{(-n)}; \boldsymbol{\beta}\boldsymbol{\alpha} \right)_{\boldsymbol{k}\boldsymbol{j}} \, \in \, \left(\sigma^{\Gamma}\right)^{-n} \left(\mathfrak{X}^{+} \left( z_{0}; \boldsymbol{\alpha} \right)_{\boldsymbol{j}}\right)} \delta_{\mathfrak{X}^{+} \left( z_{k_{n}}^{(-n)}; \boldsymbol{\beta}\boldsymbol{\alpha} \right)_{\boldsymbol{k}\boldsymbol{j}}} \ \ \rightharpoonup\ \ \mu.
    \end{equation}
    It is easy to observe that $(\Pi_{0})_{*} \mu_{n} = \dfrac{1}{d_{{\rm top}}^{n}} \displaystyle \sum \limits_{w \, \in \, (F^{n})^{\dagger} (z)} \delta_{w}$, the limit of which is the Dinh-Sibony measure $\mu_{\digamma}$. As we have $(\Pi_{0})_{*} \mu_{n} \to (\Pi_{0})_{*} \mu$, we get $(\Pi_{0})_{*} \mu = \mu_{\digamma}$. Theorem 7.1 in \cite{lop:1994} also guarantees that $h_{\mu} (\sigma^{\Gamma}) = \log d_{{\rm top}}$ and that for any $\mu' \neq \mu$ we have $h_{\mu'} (\sigma^{\Gamma}) < \log d_{{\rm top}}$, where $h_{\mu} (\sigma^{\Gamma})$ denotes the entropy of the map $\sigma^{\Gamma}$ with respect to the measure $\mu$. Now, from the proof of Theorem 1.1 in \cite{bks:2026}, we know that $(\Pi_{0})_{*}\mu$ is the unique measure of maximal entropy for the correspondence $\digamma$, which completes the proof.
\end{proof}

The moral of the story from above proof is the following. The work \cite{tmss:2026} dealt with studying ergodicity and initiating the concept of induced correspondences. Some of the main results of that work, most notably Theorem 5.1, which is actually a version of Kac's theorem in the setting of holomorphic correspondences, use the following space of measures. 
\[
\mathscr{E}^{\Gamma} = \left\{ \nu \in \mathcal{M} (\widehat{\mathbb{C}}): \nu = (\Pi_{0})_{*}\mu, \ \text{for some} \ \sigma^{\Gamma}-\text{ergodic}\ \mu \in \mathcal{M} \left(\mathscr{P}^{\Gamma} \left( \widehat{\mathbb{C}} \right), \sigma^{\Gamma} \right)\right\}.
\]
Theorem 5.1 in \cite{tmss:2026} holds for measures $\mu \in \mathscr{E}^{\Gamma}$ such that $\mu$ is also $\digamma^{*}$ invariant. Results of Dinh and Sibony (see for example \cite{ds:2006}) show that the measure $\mu_{\digamma}$ is $\digamma^{*}$ invariant. Further, by the construction followed in our proof above, $\mu_{\digamma} \in \mathscr{E}^{\Gamma}$. Thus, our result gives an interesting example of a measure in the space $\mathscr{E}^{\Gamma}$ and hence a situation where appropriate results from \cite{tmss:2026} can be applied. Further, recall from \cite{lon:2024, dkw:2020, sb:2026} that the correspondence $\digamma$ is ergodic with respect to $\mu_{\digamma}$ and it also satisfies the mixing property with respect to that measure. 

\section{Some discussion on Lemma \ref{lem:lf}}
\label{sec:seven}

We end this paper with a short discussion on the strategy used in the proof of Lemma \ref{lem:lf}. This has to do with the choice of the rate function $I_{Q_{F}}$ as seen in that proof. Note that we had to use a rate function defined on $\mathcal{M} \left( \mathscr{P}^{\Gamma}_{{\rm inv}} \left( \Omega_{{\rm DS}} \right) \right)$ and not on $\mathcal{M} \left( \Omega_{{\rm DS}} \right)$, even as the statement in Lemma \ref{lem:lf} studies terms involving points in the ambient space $\Omega_{{\rm DS}}$. The other approach of using a rate function defined on $\mathcal{M} \left( \Omega_{{\rm DS}} \right)$ has a few shortcomings with regards to our situation, which we record below. 

For a H\"{o}lder continuous function $f : \Omega_{{\rm DS}} \longrightarrow \mathbb{R}$ and define a function $Q_{f} : \mathcal{C} \left(\Omega_{{\rm DS}}, \mathbb{R} \right) \longrightarrow \mathbb{R}$ given by $Q_{f} (g) = {\rm Pr} (\Gamma, f + g) - {\rm Pr} (\Gamma, f)$. For any measure $\nu \in \mathcal{M} \left( \Omega_{{\rm DS}} \right)$ define, the real valued function $I_{Q_{f}}$ on $\mathcal{M} \left( \Omega_{{\rm DS}} \right)$ given by,
\begin{equation}
\label{eqn:iq}
    I_{Q_{f}} (\nu) := \sup\limits_{g\, \in\, \mathcal{C} \left(\Omega_{{\rm DS}}, \mathbb{R} \right)} \left( \int g\, \mathrm{d}\nu - Q(g) \right).
\end{equation}
Note that if we have $F = f \circ \Pi_{0}$ and $G = g \circ \Pi_{0}$, we actually get $Q_{F} (G) = Q_{f} (g)$. The main problem is that the function $Q_{f}$ on $\mathcal{C} \left(\Omega_{{\rm DS}}, \mathbb{R} \right)$ need not extend to the function $Q_{F}$ on $\mathcal{C} (\mathscr{P}^{\Gamma}_{{\rm inv}} \left( \Omega_{{\rm DS}} \right), \mathbb{R})$, for the choice $F = f \circ \Pi_{0}$. This happens because not all functions $\mathcal{C} (\mathscr{P}^{\Gamma}_{{\rm inv}} \left( \Omega_{{\rm DS}} \right), \mathbb{R})$ are of the form $g \circ \Pi_{0}$, for some $g \in \mathcal{C} \left(\Omega_{{\rm DS}}, \mathbb{R} \right)$.

Further, upon simplification of the expression in Equation \eqref{eqn:iq}, we observe that, 
\begin{equation}
\label{eqn:iqsimple}
    I_{Q_{f}} (\nu) = - \inf\limits_{g\, \in\, \mathcal{C} \left(\Omega_{{\rm DS}}, \mathbb{R} \right)} \left({\rm Pr} \left( \Gamma,\; g \right) - \displaystyle \int g\, \mathrm{d}\nu \right) + {\rm Pr} \left( \Gamma,\; f \right) - \displaystyle \int f\, \mathrm{d}\nu.
\end{equation}
Note that, the map $\nu \mapsto h_{\nu} (\Gamma)$ need not be upper semi-continuous (even under the assumptions that $\digamma$ is distance expanding) on $\mathscr{S}^{\Gamma}$, unlike the analogous situation of rational maps on their Julia set. This difficulty arises primarily from the very definition of $h_{\nu} (\Gamma)$ (see Definition \ref{defn:entcorr}), which involves a certain supremum, making the upper semi-continuity difficult to achieve. On this note, it is interesting to recall a result from \cite{ss:2026} which gives an equivalent criterion for the upper semi continuity of the entropy map $\nu \mapsto h_{\nu} (\Gamma)$ and is relevant in our context. The result showed that $h_{\nu'} (\Gamma) = \inf \limits_{g \, \in \, \mathcal{C} (\Omega_{{\rm DS}}, \mathbb{R})} \left \{ {\rm Pr} (\Gamma, g) - \displaystyle \int g \, \mathrm{d}\nu' \right\}$ if and only if the map $\nu \mapsto h_{\nu} (\Gamma)$ is upper semi continuous at $\nu' \in \mathscr{S}^{\Gamma}$. 

If the entropy map associated to the correspondence is upper semi-continuous, then using Equation \eqref{eqn:iqsimple}, we obtain $I_{Q_{f}}$ to be lower semi-continuous on $\mathscr{S}^{\Gamma}$. Even if $I_{Q_{f}}$ is not lower semi-continuous, one can consider a lower semi-continuous regularisation of $I_{Q_{f}}$, given by $J_{Q_{f}} (\nu) := \lim \limits_{n \, \to \, \infty} \left( \inf \left\{ I_{Q_{f}} (\nu') : \nu' \in B \left(\nu, \frac{1}{n} \right)  \right\} \right).$ But $J_{Q_{f}}$ is not directly related to the thermodynamic formalism of correspondences and hence, is not entirely interesting. So, we use the power of the bijection between the sets $\mathscr{Q}_{n}^{\Gamma} \left( z_{0} \right)$ and $\left(\sigma^{\Gamma}\right)^{-n} \left(\mathfrak{X}^{+} \left( z_{0}; \boldsymbol{\alpha} \right)_{\boldsymbol{j}}\right)$ and consider the function $I_{Q_{F}}$ to be the rate function for our problem. It also turns out as we observed in Section \ref{sec:four} that the rate function $I_{Q_{F}}$ has some good connection to the thermodynamic formalism for correspondences - yet another reason in favour of our choice.


\begin{thebibliography}{99} 

\bibitem{bs:2016}
{\sc Bharali, G.} and {\sc Sridharan, S.}, ``The dynamics of holomorphic correspondences of $\mathbb{P}^{1}$: invariant measures and the normality set", \emph{Complex Var. Elliptic Equ.} {\bf 61}, (2016), (1587 - 1613). 

\bibitem{bs:2021}
{\sc Bharali, G.} and {\sc Sridharan, S.}, ``The entropy of holomorphic correspondences: exact computations and rational semigroups", \emph{Ann. Fenn. Math.} {\bf 46}, (2021), (1135 - 1154). 

\bibitem{cr:2011}
{\sc Comman, H.} and {\sc Rivera-Letelier, J.}, ``Large deviation principles for non-uniformly hyperbolic rational maps", \emph{Ergodic Theory Dynam. Systems} {\bf 31}, (2021), (321 - 349). 

\bibitem{dkw:2020}
{\sc Dinh, T.-C., Kaufmann, L.} and {\sc Wu, H.}, ``Dynamics of holomorphic correspondences on Riemann surfaces",  \emph{Internat. J. Math.} {\bf 31}, (2020), (21 pp). 

\bibitem{ds:2006}
{\sc Dinh, T-C.} and {\sc Sibony, N.}, ``Distribution des valeurs de transformations méromorphes et applications", [Distribution of the values of meromorphic transformations and applications],  \emph{Comment. Math. Helv.}, {\bf 81} (2006), (221 - 258). 

\bibitem{ds:2008}
{\sc Dinh, T.-C.} and {\sc Sibony, N.}, ``Upper bound for the topological entropy of a meromorphic correspondence", \emph{Israel J. Math. } {\bf 163}, (2008), (29 - 44). 

\bibitem{ss:2026}
{\sc Gopinathan, S.} and {\sc Sridharan, S.}, ``A variational principle for holomorphic correspondences", \emph{Internat. J. Math.}, (To appear) (17 pp.). 

\bibitem{hem:2024}
{\sc Hemmingsson, N.}, ``Equidistribution of iterations of holomorphic correspondences and Hutchinson invariant sets", \emph{Conform. Geom. Dyn.} {\bf 28}, (2024), (97 - 114). 

\bibitem{llz:2023}
{\sc Li, X., Li, Z.,} and {\sc Zhang, Y.}, ``Thermodynamic formalism for correspondences", \texttt{arXiv:2311.09397}, (49 pp.).

\bibitem{lju:1983}
{\sc Lyubich, M.Ju.}, ``Entropy properties of rational endomorphisms of the Riemann sphere", \emph{Ergodic Theory Dynam. Systems}, {\bf 3} (1983), (351 - 385).

\bibitem{lon:2022}
{\sc Londhe, M.}, ``Recurrence in the dynamics of meromorphic correspondences and holomorphic semigroups", \emph{Indiana Univ. Math. J.}, {\bf 71} (2022), (1131 - 1154). 

\bibitem{lon:2024}
{\sc Londhe, M.}, ``Ergodicity in the dynamics of holomorphic correspondences", \emph{Ann. Fenn. Math.}, {\bf 49} (2024), (695 - 712). 

\bibitem{lop:1994}
{\sc Lopes, A.}, ``Entropy, pressure and large deviation", \emph{Cellular automata, dynamical systems and neural networks (Santiago, 1992)}, \emph{Math. Appl.}, {\bf 282} (1994), (79 - 146).

\bibitem{lop:1990}
{\sc Lopes, A.O.}, ``Entropy and large deviation", \emph{Nonlinearity}, {\bf 3} (1990), (527 - 546). 

\bibitem{sb:2026}
{\sc Mana, S.T.} and {\sc Seshadri, B.K.}, ``A view towards mixing in holomorphic correspondences", \texttt{arXiv:2604.00042}, (17 pp.). 

\bibitem{tmss:2026}
{\sc Mana, S.T.} and {\sc Sridharan, S.}, ``On induced systems and ergodicity for holomorphic correspondences", \texttt{arXiv:2609.29439}, (23 pp.). 

\bibitem{vd:2024}
{\sc Matus de la Parra, V.}, ``Equidistribution for matings of quadratic maps with the modular group", \emph{Ergodic Theory Dynam. Systems} {\bf 44}, (2024), (859 - 887). 

\bibitem{ps:1996}
{\sc Pollicott, M.} and {\sc Sharp, R.}, ``Large deviations and the distribution of preimages of rational maps", \emph{Comm. Math. Phys.} {\bf 181}, (1996), (733 - 739). 

\bibitem{ps:2007}
{\sc Pollicott, M.} and {\sc Sridharan, S.}, ``Large deviation results for periodic points of a rational map", \emph{J. Dyn. Syst. Geom. Theor.} {\bf 5}, (2007), (69 - 77). 

\bibitem{pu:2010}
{\sc Przytycki, F.} and {\sc Urbański, M.}, ``Conformal fractals: ergodic theory methods", London Math. Soc. Lecture Note Ser., \emph{Cambridge University Press, Cambridge}, (2010), (354 pp). 

\bibitem{bks:2026}
{\sc Seshadri, B.K.} and {\sc Sridharan, S.}, ``Uniqueness of the equilibrium state in the dynamics of holomorphic correspondences", \texttt{arXiv:2607.28689}, (21 pp.). 

\bibitem{ss:2025}
{\sc Sridharan, S.} and {\sc Subith, G.}, ``A Ruelle operator for holomorphic correspondences", \emph{Conform. Geom. Dyn.} {\bf 29}, (2025), (143 - 168).

\bibitem{ssa:2021}
{\sc Sridharan, S., Tikekar. S.N.} and {\sc Tiwari, A.R.}, ``Distribution of typical orbits for random dynamics generated by finitely many rational maps", \emph{Complex Anal. Oper. Theory} {\bf 15}, (2021), (22 pp).

\bibitem{wal:1982}
{\sc Walters, P.}, ``An introduction to ergodic theory", \emph{Springer-Verlag, New York-Berlin}, (1982), (250 pp).

\end{thebibliography}
\end{document}